\documentclass[aap,MSNbibl,seceqn,dvips]{arximspdf}
\usepackage{graphics}

\doi{10.1214/09-AAP641}
\volume{20}
\issue{3}
\pubyear{2010}
\firstpage{951}
\lastpage{977}

\makeatletter
\newtheorem{lm}{Lemma}[section]
\newtheorem{prop}{Proposition}
\newproclaim{rem}{Remark}[section]

\newcommand{\overld}{\overline{\mathfrak{d}}}
\newcommand{\underld}{\underline{\mathfrak{d}}}

\makeatother

\begin{document}
\begin{frontmatter}

\title{On collisions of Brownian particles\thanksref{T1}}
\runtitle{Collisions of Brownian particles}
\thankstext{T1}{Supported in part by 
NSF Grant NSF-DMS-06-01774 at
Columbia University; the second author is on leave from Columbia
University's Department of Mathematics.}

\begin{aug}
\author[A]{\fnms{Tomoyuki} \snm{Ichiba}\ead[label=e1]{ichiba@pstat.ucsb.edu}\corref{}} \and
\author[B]{\fnms{Ioannis} \snm{Karatzas}\ead[label=e2]{ik@enhanced.com}}
\runauthor{T. Ichiba and I. Karatzas}
\affiliation{University of California, Santa Barbara and INTECH}
\address[A]{Department of Statistics \& Applied Probability \\
South Hall\\
University of California, Santa Barbara \\
Santa Barbara, California 93106 \\
USA\\
\printead{e1}}
\address[B]{INTECH Investment Management \\
One Palmer Square, Suite 441 \\
Princeton, New Jersey 08542 \\
USA\\
\printead{e2}}
\end{aug}

\received{\smonth{10} \syear{2008}}
\revised{\smonth{4} \syear{2009}}

%
\begin{abstract}
We examine the behavior of $n$ Brownian particles diffusing on the real
line with bounded,
measurable drift and bounded, piecewise continuous diffusion
coefficients that depend on the
current configuration of particles. Sufficient conditions are
established for the absence
and for the presence of triple collisions among the particles. As an
application to the
Atlas model for equity markets, we study a special construction of such
systems of
diffusing particles using Brownian motions with reflection on
polyhedral domains.
\end{abstract}

%
\begin{keyword}[class=AMS]
\kwd[Primary ]{60G17}
\kwd{60G44}
\kwd[; secondary ]{60G85}.
\end{keyword}
\begin{keyword}
\kwd{Martingale problem}
\kwd{triple collision}
\kwd{effective dimension}
\kwd{Bessel process}
\kwd{reflected Brownian motion}
\kwd{comparison theorem}
\kwd{Atlas model}.
\end{keyword}
\end{frontmatter}

\section{Introduction}\label{intro}
It is well known that, with probability one, the $n$-dimension\-al
Brownian motion started away from the origin will hit the origin
infinitely often for $n=1$ while it will never hit the origin for
$n\ge2$. This is also true for the $n$-dimensional Brownian motion
with constant drift and diffusion coefficients, by Girsanov's
theorem and re-orientation of coordinates. The next step of
generalization is the case of bounded drift and diffusion
coefficients. The existence of weak solutions for the stochastic
equations that describe such
processes was discussed by Krylov \cite{Krylov80} and Stroock and Varadhan \cite{SV06} through the study of
appropriate martingale problems.

Now let us suppose that $\mathbb{R}^{n}$ is partitioned as a finite union
of disjoint polyhedra. Bass and Pardoux \cite{BP87} established
the existence and uniqueness of a weak solution to the
stochastic integral equation,
%
%
\begin{equation}\label{eq: 1}
X(t) = x_0 + \int_0^t \mu(X(s))\,d s + \int_0^t\sigma(X(s))\, d
W(s),\qquad0 \le t< \infty,
\end{equation}
with initial condition $ x_0 \in\mathbb{R}^n $ where the measurable
functions $\mu{}\dvtx{}\mathbb{R}^n \to\mathbb{R}^n$ and $\sigma
{}\dvtx{} \mathbb{R}^n
\to\mathbb{R}^{n\times
n}$ are bounded, and, moreover, $\sigma$ is everywhere nonsingular
and piecewise constant (i.e., constant on each polyhedron).
The continuous process $\{W(t), 0\le t < \infty\}$ is
an $n$-dimensional Brownian motion on some filtered
probability space $(\Omega, \mathcal{F}, \{\mathcal{F}_t\},\mathbb
{P})$. Here uniqueness is
understood in the
sense of the probability distribution.

Bass and Pardoux also discovered an interesting phenomenon,
namely, that the weak solution to (\ref{eq: 1}) may satisfy
%
%
\begin{equation}\label{eq: 2}
\mathbb{P}_{x_0} \bigl(X(t) = 0, \mbox{ i.o.}\bigr) = 1;\qquad
x_0\in\mathbb{R}^n,
\end{equation}
for a diffusion matrix $\sigma(\cdot)$ with special structure and
without drift $ \mu(\cdot) \equiv0 $.
Here $\mathbb{P}_{x_0}$ is the solution to the martingale problem
corresponding to (\ref{eq: 1}). In the Bass and Pardoux \cite{BP87}
example, the
whole space $\mathbb{R}^n$ is partitioned into a finite number of
polyhedral domains with common vertex at the origin, carefully
chosen \textit{small} apertures and $\sigma(\cdot)$ constant in
each domain. We review this example in Remark~\ref{rem: bass and pardoux}.

In the present paper we find conditions sufficient for ruling
(\ref{eq: 2}) out. More specifically, we are interested in the
case of a bounded, measurable drift vector $ \mu(\cdot) $ and
of a bounded, piecewise continuous diffusion matrix, 
%
%
\begin{equation}\label{eq: 3}
\sigma(x) = \sum_{\nu=1}^m \sigma_\nu(x)
\mathbf{1}_{\mathcal{R}_\nu}(x) \equiv\sigma_{\mathfrak{p}(x)}
(x);\qquad
x \in\mathbb{R}^n,
\end{equation}
under the assumption of \textit{well-posedness} (existence and
uniqueness of solution) when $n \ge3 $.
Here $\mathbf{1}_{\{\cdot\}}$ is the indicator function; the sets
$ \{\mathcal{R}_\nu\}_{ \nu= 1}^m $\vspace*{1pt} form a partition of $\mathbb
{R}^n$ for
some $ m \in\mathbb{N} $, namely, $\mathcal{R}_\nu\cap\mathcal
{R}_\kappa=
\varnothing$ for $ \nu\ne\kappa$ and $\bigcup_{\nu=1}^m \mathcal
{R}_\nu
= \mathbb{R}^n$
and the mapping $\mathfrak{p}{}\dvtx{} \mathbb{R}^n \to\{1, \ldots
, m\}$
satisfies $ x \in\mathcal{R}_{\mathfrak{p}(x)}$ for every $ x\in
\mathbb{R}^n $. Throughout this paper we shall assume that $\mathcal
{R}_{\nu}$
is an $n$-dimensional polyhedron for each $\nu= 1, \ldots, m $, and
that the $(n \times n)$ matrix-valued functions $ \{\sigma_\nu(\cdot)
\sigma_\nu^{\prime} (\cdot) \}_{
\nu= 1}^m $ are positive-definite everywhere.

We shall also assume throughout that there exists a unique weak
solution for equation (\ref{eq: 1}). Existence is guaranteed by the
measurability and boundedness of the functions $ \mu(\cdot) $ and $
\sigma(\cdot) \sigma' (\cdot) $ as well as the uniform strong
nondegeneracy of $ \sigma(\cdot) \sigma' (\cdot) $ (e.g., Krylov
\cite{MR2078536}, Remark~2.1) where the superscript $\prime$
represents the
transposition. Uniqueness holds
when $ n = 1$ or $n= 2 $; for $n \ge3 $, the argument of Chapter~7 of
Stroock and Varadhan \cite{SV06} implies uniqueness if the function
$\sigma(\cdot)$ in (\ref{eq: 3}) is continuous on $\mathbb{R}^{n}$
(Theorem~7.2.1
of \cite{SV06}) or close to constant (Theorem~7.1.6 of \cite{SV06}),
namely, if there exists a constant $(n\times n)$ matrix $\alpha
$ and a sufficiently small $\delta>0 $, depending on the dimension $n$
and the bounds of eigenvalues of $\sigma(\cdot)$ such that
$\sup_{x \in\mathbb{R}^{n}} \Vert\sigma(x) \sigma^{\prime} (x) -
\alpha
\Vert\le\delta$. Bass and Pardoux \cite{BP87} showed uniqueness for
piecewise-constant coefficients, that is, $\sigma_{\nu}(\cdot)
\equiv
\sigma_{\nu} $, $\nu= 1, \ldots, m$. For further discussion on
uniqueness and non-uniqueness, we refer to the paper by Krylov \cite
{MR2078536} and the references therein. The structural assumption
(\ref{eq: 3}) may be weakened to more general bounded cases, under
modified conditions.

Our main concern is to obtain sufficient conditions on
$\mu(\cdot)$ and on $\sigma(\cdot)$ of the form (\ref{eq: 3}) so
that with $ n \ge3 $ we have
\begin{eqnarray}\label{eq: 4}
\mathbb{P}_{x_0} \bigl(X_i(t) = X_j(t) = X_k(t), \mbox{ for some }
t\ge0
\bigr) &=& 0\quad\mbox{or}\nonumber\\[-8pt]\\[-8pt]
\mathbb{P}_{x_0} \bigl(X_i(t) = X_j(t) = X_k(t), \mbox{ for some } t
\ge0
\bigr) &=&1;\qquad x_{0}\nonumber
\in\mathbb{R}^{n},
\end{eqnarray}
for some $1\le i< j < k \le n $. Put differently, we study
conditions on their drift and diffusion coefficients, under
which three Brownian particles moving on the real line can collide
at the same time, and conditions under which such ``\textit{triple
collisions}'' can never occur. Propositions~\ref{prop1} and~\ref
{prop2} provide
answers to these questions in Section~\ref{sec: first approach}.

In Section~\ref{sec: second approach} we study a class of the weak
solutions to the stochastic differential equation (\ref{eq: 1}),
clarifying the relationship between the rank of process coordinates and
the reflected Brownian motion on $ (n-1)$-dimensional polyhedral
domain. Proposition~\ref{prop3} shows that the process has no triple collisions
under some parametric conditions.

The results have consequences in the computations of local
times for the differences $\{X_i(t)-X_j(t), X_j(t)-X_k(t)\}$. We
discuss such local times with application to the analysis of a
so-called ``Atlas model'' for equity markets in Section~\ref{sec: application}.
Proofs of selected results are presented in \hyperref[sec: appendix]{Appendix}.

Recent work related to this problem was done by C\'{e}pa and
L\'{e}pingle \cite{CL07}. These authors consider a system of
mutually repelling Brownian particles and show the absence of
triple collisions. The electrostatic repulsion they consider comes from
unbounded drift coefficients; in our setting, all drifts are
bounded.

\section{A first approach}\label{sec: first approach}
\subsection{The setting}

Consider the stochastic integral equation (\ref{eq: 1}) with coefficients
$\mu(\cdot)$ and $\sigma(\cdot)$ as in (\ref{eq: 3}), and assume
that the matrix-valued functions $\sigma_\nu(\cdot)$, $\nu= 1,
\ldots,
m$, are uniformly positive-definite. Then the inverse $\sigma
^{-1}(\cdot
)$ of the
diffusion coefficient $\sigma(\cdot)$ exists in the sense $
\sigma^{-1}(\cdot) = \sum_{\nu=1}^m \sigma_\nu^{-1} (\cdot)
\mathbf
{1}_{\mathcal{R}_\nu} (\cdot) $. As usual, a~\textit{weak solution}
of this
equation consists of a probability space $ (\Omega, \mathcal{F},
\mathbb
{P}) $; a filtration $\{\mathcal{F}_t, 0 \le t < \infty\}$ of
sub-$\sigma$-fields of $\mathcal{F}$ which satisfies the ``usual conditions'' of
right-continuity and augmentation by the
$\mathbb{P}$-negligible sets in $\mathcal{F}$; and two adapted, $n$-dimensional
processes on this space $ X(\cdot), W(\cdot) $ on this space, such
that $ W(\cdot) $ is Brownian motion and (\ref{eq: 1}) is satisfied $
\mathbb{P}$-almost surely. The concept of uniqueness associated with this
notion of solvability, is \textit{uniqueness in distribution} for
$X(\cdot) $.

\subsection{Removal of drift}\label{sec: removal of drift}

We start by observing that the bounded drift
has no effect on the probability of absence of triple collisions.
Indeed, if we define an $n$-dimensional process
$
\xi(t):= \sigma^{-1}(X(t))\mu(X(t)), 0\le t < \infty$,
then the nature of the functions $\mu(\cdot)$ and $\sigma(\cdot)$ in
(\ref{eq: 3}) guarantees that the mapping $t \mapsto\xi(t)$ is
right-continuous
or left-continuous on each boundary $ \partial
\mathcal{R}_{\mathfrak{p}(X(t))} $ at each time $t$, deterministically,
according to the position $\mathcal{R}_{\mathfrak{p}(X(t-))}$ of
$X(t-)$. Thus,
although the sample path of $n$-dimensional process $\xi(\cdot)$ is
not entirely right-continuous or left-continuous, it \textit{is}
progressively measurable. Moreover, $\xi(\cdot)$ is bounded, so
the exponential process
%
%
\begin{equation}\label{eq: eta}
\hspace*{20pt}\eta(t) = \exp\biggl[ - \int_0^t \langle\xi(u),
d W (u) \rangle- \frac{1}{2}\int^t_0 \Vert\xi(u)\Vert^2\,d u
\biggr];\qquad0 \le t < \infty,
\end{equation}
is a continuous martingale where $\Vert x \Vert^2:=
\sum_{j=1}^n x_j^2, x \in\mathbb{R}^n$, stands for $n$-dimensional
Euclidean norm, and the bracket $ \langle x, y \rangle:=
\sum_{j=1}^{n} x_{j} y_{j} $ is the inner product of two vectors
$x, y \in\mathbb{R}^{n}$. By Girsanov's theorem,
%
%
\begin{equation}\label{eq: BM from Grisanov}
\widetilde W(t):= W(t) + \int^t_0 \sigma^{-1}(X(u))\mu(X(u))\, d
u,\qquad
\mathcal{F}_t; 0 \le t < \infty,
\end{equation}
is an $n$-dimensional standard Brownian motion under the
new probability measure $\mathbb{Q}$, \textit{locally} equivalent to $
\mathbb
{P} $, that satisfies
%
%
\begin{equation}\label{eq: def of QQ}
\mathbb{Q}_{x_{0}} (C) = \mathbb{E}^{\mathbb{P}_{x_{0}}} (\eta
(T)1_C );\qquad C\in\mathcal{F}_T,
0\le T <
\infty.
\end{equation}

Let us define an increasing family of events $ C_{T}:= \{ X_{i}(t) =
X_{j}(t) = X_{k}(t)$, for some $t \in[0, T] \} $, $ T \ge0 $.
If we know a priori that
%
%
\begin{equation}\label{eq: no hit under Q}
\mathbb{Q}_{x_{0}}\bigl(X_{i}(t) = X_{j}(t) = X_{k}(t), \mbox{ for
some } t
\ge
0 \bigr) = 0,
\end{equation}
then we obtain $ 0 = \mathbb{Q}_{x_{0}}(C_{\ell}) = \mathbb
{P}_{x_{0}}(C_{\ell}) $
for $ \ell\ge1 $, and so
\begin{eqnarray}\label{eq: no hit under P}
\mathbb{P}_{x_{0}}\bigl(X_{i}(t) = X_{j}(t) = X_{k}(t), \mbox{ for
some } t
\ge0\bigr)
&=& \mathbb{P}_{x_{0}}\Biggl(\bigcup_{\ell=1}^{\infty} C_{\ell
}\Biggr)
\nonumber\\[-8pt]\\[-8pt]
&=&\lim_{\ell\to\infty} \mathbb{P}_{x_{0}}(C_{\ell}) =
0.\nonumber
\end{eqnarray}

Thus, in order to evaluate the probability of absence of triple
collisions in (\ref{eq: 4}), it is enough to consider the case of $\mu
(\cdot) \equiv0$ in
(\ref{eq: 1}), namely
%
%
\begin{equation}\label{eq: 5}
X(t) = x_0 + \int_0^t \sigma(X(s) ) \,d \widetilde{W} (s),\qquad0
\le t <\infty,
\end{equation}
under the new probability measure $ \mathbb{Q}_{x_{0}} $.
The infinitesimal generator $\mathcal A$ of this process, defined
on the space $C^{2}(\mathbb{R}^{n}; \mathbb{R})$ of twice continuously
differentiable functions $ \varphi{}\dvtx{} \mathbb{R}^{n}\to
\mathbb{R}$, is
given as
%
%
\begin{equation}\label{eq: 6}
\mathcal A \varphi(x):= \frac{1}{2}\sum_{i,k = 1}^n a_{ik}(x) \frac
{\partial^2}{\partial x_i\, \partial x_k} [\varphi(x) ];\qquad
\varphi\in C^{2}(\mathbb{R}^{n}; \mathbb{R}),
\end{equation}
where $\sigma_{ij}(\cdot)$ is the $(i,j)$th element of
the matrix-valued function $\sigma(\cdot)$, and
%
%
\begin{equation}\label{eq: 6aA}
\qquad a_{ik}(x):= \sum_{j = 1}^n \sigma_{ij}(x)\sigma_{kj}(x),\qquad
A(x):= \{a_{ij}(x) \}_{1\le i,j\le n};\qquad
x \in\mathbb{R}^{n}.
\end{equation}
The uniform positive-definiteness
of the matrices $ \{\sigma_\nu\sigma_\nu^{\prime}\}(\cdot), \nu= 1,
\ldots, m $, in (\ref{eq: 3}) implies that the operator $ \mathcal A $
is uniformly elliptic. As is well known from \cite{SV06}, existence
(respectively, uniqueness) of a weak solution to the stochastic
integral equation (\ref{eq: 5}), is equivalent to the solvability
(respectively, well-posedness) of the martingale problem associated
with the operator $ \mathcal{A} $.

\subsection{Comparison with Bessel processes}\label{sec: equivalent problems}

Without loss of generality we start from the case
$i = 1, j = 2, k = 3$ in (\ref{eq: 4}). Let us define $(n\times1)$
vectors $d_1, d_2, d_3$ to extract the information
of the diffusion matrix $\sigma(\cdot)$ on $(X_{1}, X_{2}, X_{3})$,
namely
\begin{eqnarray*}
d_1&:=& (1, -1, 0,\ldots, 0)^\prime,\qquad
d_2:= (0, 1, -1, 0,\ldots,0)^\prime,
\\
d_3&:=& (-1, 0, 1, 0,\ldots,0)^\prime,
\end{eqnarray*}
where the superscript $ \prime$ stands for transposition.
Define the $(n \times3)$-matrix $D = (d_1, d_2, d_3)$ for notational
simplicity. The
cases we consider in (\ref{eq: 4}) for $ i = 1, j = 2, k =
3 $ are equivalent to
\begin{eqnarray*}
\mathbb{P}_{x_0} \bigl( s^{2}(X(t)) = 0, \mbox{ for some } t \ge0
\bigr)
&=& 0\quad\mbox{and }\nonumber\\[-8pt]\\[-8pt]
\mathbb{P}_{x_0} \bigl( s^{2}(X(t)) = 0, \mbox{ for some } t \ge0
\bigr)
&=& 1;\qquad
x_{0}\in\mathbb{R}^{n}\nonumber
\end{eqnarray*}
where the continuous function
\begin{eqnarray}\label{eq: 7s}
s^2(x) &:=& (x_1-x_2)^2+(x_2 - x_3)^2 + (x_3 - x_1)^2 \nonumber\\
[-8pt]\\[-8pt]
&=& d_1^\prime x x^\prime d_1 + d_2^\prime x x^\prime d_2 +
d_3^\prime x x^\prime d_3 = x^\prime D D^\prime x;\qquad x\in\mathbb{R}
^n\nonumber
\end{eqnarray}
measures the sum of squared distances for the three
particles of interest. Thus, it suffices to study the
behavior of the continuous, nonnegative process $\{s^2(X(t));
0 \le t< \infty\}$ around its zero set
%
%
\begin{equation}\label{eq: 9}
\mathcal Z:= \{ x \in\mathbb{R}^{n}{}\dvtx{} s(x) = 0\}.
\end{equation}

Let us define the following positive, piecewise continuous functions $
Q(\cdot) $, $ \widetilde{R} (\cdot) $ computed from the
variance--covariance matrix $ A(\cdot) = \sigma(\cdot) \sigma
^{\prime}
(\cdot) $:
\begin{eqnarray}\label{eq: def of Q}
\hspace*{30pt} \widetilde{R} (x)&:=& \frac{\operatorname{trace}(D^{\prime}A(x) D) \cdot
x^{\prime} D D^{\prime} x}{x^{\prime} D D^{\prime} A(x) D D^{\prime
} x}
= \frac{\operatorname{trace}(D^{\prime} A(x) D)}{Q(x)},\qquad\mbox
{where}\nonumber\\[-8pt]\\[-8pt]
Q(x)&:=& \frac{x^{\prime} D D^{\prime} A(x) D D^{\prime
}x}{x^{\prime} D D^{\prime} x};
\qquad x \in\mathbb{R}^{n} \setminus\mathcal Z.\nonumber
\end{eqnarray}

Under the new probability measure $ \mathbb{Q}_{x_{0}} $ of (\ref{eq:
def of
QQ}) the process $ s(X(\cdot)) $ is a semimartingale with decomposition
$ d s(X(t)) = \widetilde{h} (X(t))\,d t + d \widetilde{\Theta} (t) $ where
%
%
\begin{eqnarray}\label{eq: h Theta}
\widetilde{h} (x) &:=&{ 1 \over2 s^3 (x) } \Biggl(s^2(x) \sum_{i=1}^3
d_i^\prime\sigma(x) \sigma(x)^\prime d_i - \Bigg\|\sum_{i=1}^3
\sigma(x)^\prime d_i d_i^\prime x \Bigg\|^2 \Biggr)\nonumber
\\
&=& \frac{x^\prime D D^\prime x \cdot\operatorname{trace} ( D^\prime A(x) D ) - x^\prime D D^\prime A(x)
D D^\prime x}{2 (x^\prime D D^\prime x)^{3/2} }\\
& =& \frac{ (\widetilde{R} (x) - 1 )Q(x)}{2 s(x)};\qquad
x \in\mathbb{R}^{n} \setminus\mathcal Z,\nonumber
\end{eqnarray}
and
\begin{eqnarray*}
\widetilde\Theta(t)&:=& \int^t_0 \Biggl( \sum_{i=1}^3 \frac{
\sigma^\prime(X(\tau)) d_i d_i^\prime X(\tau)}{s(X(\tau))} \Biggr
)\, d\widetilde{W}(\tau), \\
\langle\widetilde\Theta\rangle(t) &=& \int^t_0 \frac{x^\prime
DD^\prime A(x) D D^\prime x}{ x^\prime D D^\prime x} \bigg|_{x=X(\tau)}
\ d \tau= \int^{t}_{0} Q(X(\tau))\, d \tau;\qquad0\le t < \infty,
\end{eqnarray*}
respectively.
Here, as we shall see (\ref{eq: lower bnd for R}) in Remark~\ref{rem: 6},
we have
%
%
\begin{equation}\label{ineq: Q}
Q(\cdot) \ge c_{0}:= 3 \min_{1\le i \le n, x \in\mathbb{R}^{n}
\setminus
\mathcal Z} \lambda_{i}(x) > 0\qquad\mbox{in } \mathbb{R}^{n}
\setminus
\mathcal Z
\end{equation}
for the eigenvalues $ \{\lambda_{i}(\cdot), 1\le i \le n\} $ of $
A(\cdot) $, and so $ \langle\widetilde\Theta\rangle(\cdot) $ is
strictly increasing, when $ X(\cdot) \in\mathbb{R}^{n} \setminus
\mathcal Z $.
Now we define the increasing family of stopping times $ \Lambda_u:=
\inf
\{ t \ge0{}\dvtx{} \langle\widetilde\Theta\rangle(t)\ge u\} $, $
0 \le u <
\infty$, and note that we have
\[
\mathfrak{s} (u):= s(X(\Lambda_u)) = s(x_0) + \int^{\Lambda_u}_0
\widetilde{h} (X(t))\,d t + \widetilde B(u);\qquad
0 \le u < \infty,
\]
where $ \widetilde B(u):= \widetilde\Theta(\Lambda_u), 0 \le u <
\infty$, is a
standard Brownian motion, by the Dambis--Dubins--Schwarz theorem on
time-change for martingales. Thus, with $ \mathfrak{d}(u):= \widetilde
{R}(X(\Lambda_u)) $ we can write
%
%
\begin{equation}\label{bessel}
d \mathfrak{s}(u) = \frac{ \mathfrak{d} (u) - 1 } {2 \mathfrak
{s}(u) }\,d u+ d \widetilde B(u);\qquad0 \le u <
\infty,
\end{equation}
because
\begin{eqnarray*}
\widetilde{h} (X(\Lambda_u) ) \Lambda'_u = \frac{ [
\widetilde{R} (X(\Lambda_u)) - 1] Q(X(\Lambda_u)) } {2 s(X(\Lambda_u))
}\cdot{ 1 \over Q(X(\Lambda_u)) } = \frac{
\mathfrak{d} (u) - 1 } {2 \mathfrak{s}(u) }.
\end{eqnarray*}
The dynamics of the process $ \mathfrak{s} (\cdot)$ are therefore
comparable to those of the $\delta$-dimensional Bessel process, namely
\[
d \mathfrak{r}(u) = \frac{ \delta-1 }{2 \mathfrak{r}(u)}\, d u +
d \widetilde B (u);\qquad0 \le u < \infty.
\]
By a comparison argument similar to Ikeda and Watanabe \cite{IW77} and
Exercise 5.2.19 in \cite{KS91}, we prove in Section~\ref{sec: pf of
comparison lem} the following result.
%
\begin{lm}\label{lm: comparison}
Suppose $ x_{0} \in\mathbb{R}^{n} \setminus\mathcal Z $.
If\ $\underld:= \operatorname{essinf}
\inf_{0 \le t<\infty} \mathfrak{d}(t) \ge2 $, 
%
%
\begin{equation}\label{eq: lem of hit 1}
\mathbb{Q}_{x_{0}} \bigl( \mathfrak{s}(t) > 0, \mbox{ \textup{for
some} }
t \ge0 \bigr) = 0.
\end{equation}
If, on the other hand, $ \overld:= \operatorname{essup}
\sup_{0 \le t
<\infty} \mathfrak{d}(t) < 2 $, then
%
%
\begin{equation}\label{eq: lem of hit 2}
\mathbb{Q}_{x_{0}} \bigl( \mathfrak{s}(t) = 0, \mbox{ \textup{for
infinitely many} } t \ge0\bigr) = 1;
\end{equation}
and we have the following estimate:
%
%
\begin{equation}\label{eq: lem of hit 3}
\mathbb{Q}_{x_{0}} \bigl( \mathfrak{s}(t) = 0, \mbox{ \textup{for
some} }
t \in[0, T] \bigr) \ge1
- \kappa(T; s(x_{0}), \overld),
\end{equation}
where $ \kappa(\cdot; y, \delta) $ is the tail distribution of the
first hitting-time at the origin for Bessel process in dimension $
\delta\in(0, 2) $, starting at $ y >0$,
%
%
\begin{equation}\label{eq: kappa}
\hspace*{20pt}\kappa(T; y, \delta):= \int^{\infty}_{T} \frac{1}{t
\Gamma( \delta)}
\biggl( \frac{y^{2}}{2 t} \biggr)^{\delta} e^{- y^{2}/2 t}\, d
t;\qquad
0 \le T <\infty, y > 0.
\end{equation}
\end{lm}

This function decreases as $ T^{-\delta} $ with $ T\uparrow\infty$.
Combining Lemma~\ref{lm: comparison} with the reasoning in
Section~\ref{sec: removal of drift} and the definition $ \mathfrak
{d}(\cdot)
= \widetilde{R} (X(\Lambda_{\cdot})) $,
we obtain the following result on the \textit{absence} of triple collisions:
%
\begin{prop}\label{prop1}
Suppose that the matrices
$ \sigma_\nu(\cdot), \nu= 1,\ldots,m$, in $(\ref{eq: 3})$ are
uniformly bounded and positive-definite and satisfy the following
condition:
%
%
\begin{equation}\label{eq: condition ED and R}
\inf_{x \in\mathbb{R}^{n} \setminus\mathcal Z} \widetilde{R} (x)
\ge2 %
\end{equation}
for $ \widetilde{R} (\cdot) $ in $(\ref{eq: def of Q}) $. Then for the
weak solution $X(\cdot)$ to $(\ref{eq: 5})$ we have
\begin{eqnarray*}\label{q1}
\mathbb{Q}_{x_0} \bigl( X_1(t) = X_2(t) = X_3(t), \mbox{ \textup
{for some}
} t \ge
0 \bigr) = 0\qquad \forall x_0 \in\mathbb{R}^n\setminus\mathcal Z.
\end{eqnarray*}
Reasoning as in (\ref{eq: no hit under Q})--(\ref{eq: no hit under
P}) for the weak solution $ X(\cdot) $ to $(\ref{eq: 1})$, we get
%
%
\begin{equation}\label{p1}
\qquad\mathbb{P}_{x_0} \bigl( X_1(t) = X_2(t) = X_3(t), \mbox{ \textup
{for some}
} t \ge0
\bigr) = 0\qquad\forall x_0 \in\mathbb{R}^n\setminus\mathcal Z.
\end{equation}
\end{prop}

A class of examples satisfying (\ref{eq: condition ED and R}) is given
in Remarks~\ref{rem: R is zero}--\ref{rem: example for prop 1} and
Section~\ref{sec: proof of example 1} below. On the other hand,
regarding the \textit{presence} of triple collisions, we have the
following result; its proof is in Section~\ref{sec: pf of Prop 2}.
%
\begin{prop}\label{prop2}
Suppose that the
matrices $ \sigma_\nu(\cdot), \nu= 1,\ldots,m$, in $(\ref{eq:
3})$ are
uniformly bounded and positive-definite, and
%
%
\begin{equation}\label{eq: condition on R 2}
\delta_{0}:=
\sup_{x \in\mathbb{R}^{n} \setminus\mathcal Z } \widetilde{R} (x)
< 2.
\end{equation}
%
Then the weak solution $X(\cdot)$ to
$(\ref{eq: 5})$ starting at any $ x_0 \in\mathbb{R}^n $ satisfies
\begin{eqnarray*}\label{q4}
\mathbb{Q}_{x_0} \bigl(X_1(t) = X_2(t) = X_3(t), \mbox{ \textup{for some}
} t \ge0 \bigr)
= 1,
\end{eqnarray*}
and we have an estimate similar to $(\ref{eq: lem of hit 3})$,
\begin{eqnarray}\label{eq: lem of hit 5}
&&\mathbb{Q}_{x_{0}} \bigl( X_1(t) = X_2(t) = X_3(t), \mbox{ \textup{for
some} } t \in[0, T] \bigr) \nonumber\\[-8pt]\\[-8pt]
&&\qquad\ge1 - \kappa(c_{0} T; s(x_{0}), \delta_{0}).\nonumber
\end{eqnarray}
Here the distance function $ s(\cdot) $ and the tail probability $
\kappa( \cdot; y, \delta_{0}) $ are given by (\ref{eq: 7s}) and
(\ref
{eq: kappa}), now with dimension $ \delta_{0} \in(0,2) $ as in (\ref
{eq: condition on R 2}), and the positive constant $ c_{0} $ is given
by (\ref{ineq: Q}).

Moreover, if $ \delta_*:=\sup_{x \in\mathbb{R}^{n} \setminus
\mathcal Z} R(x)
< 2 $ holds for the modification
\begin{eqnarray}\label{eq: modified R}
R (x) &:=& \frac{ [ \operatorname{trace}(D^{\prime} A(x) D) + 2
x^{\prime} D
D^{\prime} \mu(x) ] \cdot x^{\prime} D D^{\prime} x}{x^{\prime} D
D^{\prime} A(x) D D^{\prime} x }\nonumber\\[-8pt]\\[-8pt]
&=& \widetilde R(x) + \frac{2 x^{\prime} D D^{\prime} \mu
(x)}{Q(x)};\qquad
x\in\mathbb{R}^{n} \setminus\mathcal Z,\nonumber
\end{eqnarray}
of the function $ \widetilde{R} (\cdot) $ in (\ref{eq: def of Q}), then
%
%
\begin{equation}\label{eq: triple collision prob one}
\mathbb{P}_{x_0} \bigl(X_1(t) = X_2(t) = X_3(t), \mbox{ \textup{for some}
} t \ge0 \bigr)
= 1,
\end{equation}
and we have an estimate similar to (\ref{eq: lem of hit 3}), (\ref{eq:
lem of hit 5}),
\begin{eqnarray}\label{eq: lem of hit 4}
&&\mathbb{P}_{x_{0}} \bigl( X_1(t) = X_2(t) = X_3(t), \mbox{ \textup{for
some} } t \in[0, T]\bigr)\nonumber\\[-8pt]\\[-8pt]
&&\qquad\ge1 - \kappa(c_{0}T; s(x_{0}), \delta_{*}).\nonumber
\end{eqnarray}
\end{prop}
%
%
\begin{rem}\label{rem: 6}
Since $A(\cdot)$ is positive-definite and
$\operatorname{rank} (D) = 2$, the matrix
$D^{\prime} A(\cdot) D$ is nonnegative-definite and
the number of its nonzero eigenvalues is equal to
$ \operatorname{rank}(D^{\prime} A(\cdot) D)=2 $. This implies
\begin{eqnarray*}\label{eq: lower bound for R}
\widetilde{R} (x) \ge\frac{\sum_{i=1}^{3}
\lambda^{D}_{i}(x)}
{\max_{1\le i \le3}
\lambda^{D}_{i}(x)} > 1;\qquad x \in
\mathbb{R}^{n} \setminus\mathcal Z,
\end{eqnarray*}
where $\{\lambda^{D}_{i}(\cdot), i = 1, 2, 3\}$
are the eigenvalues of the $(3 \times3)$ matrix\vspace*{1pt}
$D^{\prime} A(\cdot) D$. On the other hand, an upper bound for
$\widetilde{R} (\cdot)$ is given by
%
%
\begin{equation}\label{eq: upper bound for R}
\widetilde{R} (x) \le\frac{ \operatorname{trace}
(D^{\prime} A(x) D)}{3\min_{1 \le i \le n}
\lambda_{i} (x)};\qquad
x \in\mathbb{R}^{n} \setminus\mathcal Z,
\end{equation}
where $\{\lambda_{i}(\cdot), 1\le i \le n\}$ are
the eigenvalues of $A(\cdot)$. In fact,
we can verify
$D D^{\prime} D D^{\prime} = 3 D D^{\prime}$,
$\{x \in\mathbb{R}^{n}{}\dvtx{} DD^{\prime} x = 0\} = \mathcal Z$,
and so if $D D^{\prime} x \neq0 \in\mathbb{R}^{n}$, we obtain the
upper bound
(\ref{eq: upper bound for R}) for $ \widetilde{R} (\cdot)$ from
%
%
\begin{equation}\label{eq: lower bnd for R}
\min_{1\le i \le n} \lambda_{i} (x) \le
\frac{x^{\prime} D D^{\prime} A(x) D D^{\prime} x}
{x^{\prime} D D^{\prime} D D^{\prime} x}
= \frac{Q(x)}{3} =
\frac{\operatorname{trace} (D^{\prime} A(x) D )}{3\widetilde{R} (x)}.
\end{equation}
\end{rem}
%
%
\begin{rem}\label{rem: R is zero}
For the standard, $n$-dimensional Brownian motion, that is,
$\sigma(\cdot) \equiv I_n$, $n \ge3 $, the quantity $ \widetilde{R}
(\cdot) $ of (\ref{eq: def of Q}) is computed easily; $ \widetilde{R}
(\cdot) \equiv2 $.
More generally, suppose that the variance covariance rate
$A(\cdot)$ is
\begin{eqnarray*}
A(x):= \sum_{\nu=1}^{m} \bigl( \alpha_{\nu} I_{n} + \beta_{\nu}
D D^{\prime}
+ \mathbb{I} \mathbb{I}^{\prime} \operatorname{diag}(\gamma_{\nu}) \bigr) \cdot\mathbf{1}_{\mathcal{R}_{\nu
}}(x);\qquad
x \in\mathbb{R}^{n},
\end{eqnarray*}
for some scalar constants $\alpha_{\nu}$, $\beta_{\nu}$ and
$(n\times
1)$ constant vectors $\gamma_{\nu} $, $\nu= 1, \ldots, m$. Here
$\operatorname{diag}(x)$ is the $(n\times n)$ diagonal
matrix whose diagonal entries are the elements of $x \in\mathbb{R}^{n}$,
and $ \mathbb{I} $ is the $(n\times1)$ vector with all entries equal
to one.
Then $ \widetilde{R} (\cdot) \equiv2 $ in $ \mathbb{R}^{n} \setminus
\mathcal Z $
because $ \mathbb{I}^{\prime} D = (0, 0, 0)
\in\mathbb{R}^{1\times3}$ and
\[
D D^{\prime} = \frac{1}{3} DD^{\prime} D D^{\prime}
=
\pmatrix{
2 & -1 & -1 \cr
-1 & 2 & -1 & 0\cr
-1 & -1 & 2 \cr
& 0& 0 }
\in\mathbb{R}^{n\times n}.
\]
Hence, if the coefficients $\alpha_{\nu}, \beta_{\nu}$ and
$\gamma_{\nu}$, $\nu= 1, \ldots, m$, are chosen above so that
$A(\cdot)$ is positive-definite, we have (\ref{p1}).
\end{rem}

\begin{figure}

\includegraphics{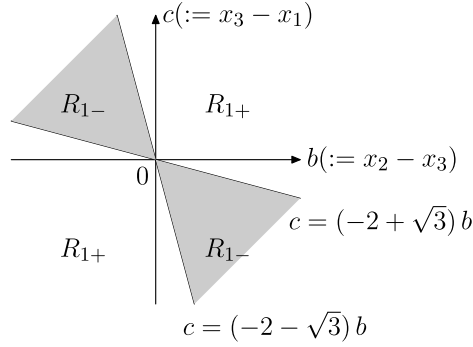}

\caption{Polyhedral regions in Remark~\protect\ref{rem: example for
prop 1}.} \label{fig: 2}
\vspace*{-5pt}
\end{figure}

%
\begin{rem}\label{rem: example for prop 1}
The condition (\ref{eq: condition ED and R}) in Proposition~\ref{prop1} holds
under several circumstances. For example,
take $n = 3$ and fix the elements
$
a_{11}(\cdot)=a_{22}(\cdot) = a_{33}(\cdot) \equiv1
$
of the symmetric matrix $A(\cdot) = \sigma\sigma^{\prime} (\cdot)$ in
(\ref{eq: 6aA}) and choose the other parameters by
%
%
\begin{eqnarray}\label{eq: def of a12}
a_{12}(x) = a_{21}(x) &:=& \alpha_{1+} \mathbf{1}_{\mathcal R_{1+}} (x)
+ \alpha_{1-} \mathbf{1}_{\mathcal R_{1-}}(x), \nonumber\\
a_{23}(x) = a_{32}(x) &:=& \alpha_{2+}\mathbf{1}_{\mathcal R_{2+}}(x)
+\alpha_{2-}\mathbf{1}_{\mathcal R_{2-}}(x), \\
a_{31}(x) = a_{13}(x) &:=& \alpha_{3+}\mathbf{1}_{\mathcal R_{3+}}(x)
+\alpha_{3-}\mathbf{1}_{\mathcal R_{3-}}(x);\qquad x \in\mathbb{R}
^{3},\nonumber
\end{eqnarray}
where $ \mathcal R_{i\pm}, i = 1, 2, 3$, are subsets of $ \mathbb
{R}^{3} $
defined by
\begin{eqnarray*}
\mathcal R_{1+} &:=& \{ x \in\mathbb{R}^{3}{}\dvtx{} {\mathfrak
f_{1}} (x) >
0 \},\qquad
\mathcal R_{2+}:= \{ x \in\mathbb{R}^{3}{}\dvtx{} {\mathfrak
f_{1}}(x) =
0,{\mathfrak f_{2}}(x) > 0 \}, \\
\mathcal R_{1-} &:=& \{ x \in\mathbb{R}^{3}{}\dvtx{} {\mathfrak
f_{1}} (x) <
0 \},\qquad
\mathcal R_{2-}:= \{ x \in\mathbb{R}^{3}{}\dvtx{} {\mathfrak
f_{1}}(x) = 0,
{\mathfrak f_{2}}(x) <0 \}, \\
\mathcal R_{3+} &:=& \{ x \in\mathbb{R}^{3}{}\dvtx{} {\mathfrak
f_{1}}(x) =
{\mathfrak f_{2}}(x) = 0, {\mathfrak f_{3}}(x) >0 \}, \\
\mathcal R_{3-} &:=& \{ x \in\mathbb{R}^{3}{}\dvtx{} {\mathfrak
f_{1}}(x) =
{\mathfrak f_{2}}(x) = 0, {\mathfrak f_{3}}(x) <0 \}, \\
{\mathfrak f_{1}} (x) &:=& \bigl[ x_{3} - x_{1} - \bigl(-2 + \sqrt
{3}\bigr) (x_{2} -x_{3})\bigr]
\cdot\bigl[ x_{3} - x_{1} - \bigl(-2 - \sqrt{3}\bigr)(x_{2 } -
x_{3})\bigr], \\
{\mathfrak f_{2}} (x) &:=& \bigl[ x_{2} - x_{3} - \bigl(-2 + \sqrt
{3}\bigr) (x_{1} -x_{2})\bigr]
\cdot\bigl[ x_{2} - x_{3} - \bigl(-2 - \sqrt{3}\bigr)(x_{1 } -
x_{2})\bigr], \\
{\mathfrak f_{3}} (x) &:=& \bigl[ x_{1} - x_{2} - \bigl(-2 + \sqrt
{3}\bigr) (x_{3} -x_{1})\bigr]
\cdot\bigl[ x_{1} - x_{2} - \bigl(-2 - \sqrt{3}\bigr)(x_{3 } -
x_{1})\bigr]
\end{eqnarray*}
for $ x \in\mathbb{R}^{3} $ with the six constants $\alpha_{i \pm} $
satisfying $ 0 < \alpha_{i+} \le1/2 $, $ -1/2 \le\alpha_{i-} < 0 $,
for $ i = 1, 2, 3 $.
Note that the zero set $\mathcal Z$ defined in (\ref{eq: 9}) is $\{
x\in\mathbb{R}^{3}{}\dvtx{} \mathfrak f_{1}(x) = \mathfrak f_{2}(x) =
\mathfrak
f_{3}(x) =0 \} $. Thus we split the region $\mathbb{R}^{3} \setminus
\mathcal
Z$ into six disjoint polyhedral regions $\mathcal R_{i\pm}, i = 1, 2, 3
$. See Figure~\ref{fig: 2}, and Section~\ref{sec: proof of example
1} for the details of this example.
\end{rem}

%

%
\begin{rem}\label{rem: bass and pardoux}
In the example of Bass and Pardoux \cite{BP87}, mentioned briefly in the
\hyperref[intro]{Introduction}, the
diffusion matrix $ \sigma(\cdot) = \sum_{\nu=1}^{m} \sigma_{\nu}
(\cdot
) \mathbf{1}_{\mathcal{R}_{\nu}} (\cdot) $
in (\ref{eq: 3}) has a special characteristic in
the allocation of its eigenvalues: {\sl All eigenvalues but the
largest are small}; namely, they are of the form $(1, \varepsilon
,\ldots,
\varepsilon)$ where $ 0 < \varepsilon< 1/2 $ satisfies, for some $ 0 <
\delta< 1/2 $,
\begin{eqnarray*}%
\bigg| \frac{x^\prime\sigma(x) \sigma^\prime(x) x} {\Vert x\Vert
^2} - 1 \bigg|
\le\delta\qquad\mbox{for } x \in
\mathbb{R}^{n} \setminus\{0\}\quad\mbox{and}\quad\frac
{(n-1)\varepsilon
^2+\delta}{1-\delta} < 1.
\end{eqnarray*}
This is the case when the diffusion matrix $ \sigma(\cdot) $ can be
written as
a piecewise constant function $ \sum_{\nu= 1} ^{m} \sigma_{\nu}
\mathbf{1}_{\mathcal{R}_{\nu}}(\cdot) $ where the constant $
(n\times n) $ matrices
$ \{\sigma_{\nu}, \nu= 1, \ldots, m \} $ have the decomposition,
\[
\sigma_{\nu} \sigma_{\nu}^{\prime}:= ( y_{\nu}, B_{\nu} )
\operatorname{diag}
(1, \varepsilon^{2}, \ldots, \varepsilon^{2} )
\pmatrix{
y_{\nu}^{\prime} \vspace*{3pt}\cr
B_{\nu}^{\prime}
},
\]
the fixed $(n\times1)$ vector $ y_{\nu}\in\mathbb{R}_{\nu} $ satisfies
\begin{eqnarray*}
\Vert y_{\nu} \Vert= 1,\qquad
\frac{\vert\langle x, y_{\nu} \rangle\vert^{2}}{\Vert x \Vert
^{2}} \ge1 - \varepsilon;\qquad
x \in\mathcal{R}_{\nu}\setminus\{0\},
\end{eqnarray*}
and the $ ( n \times(n-1)) $ matrix $ B_{\nu} $ consists of $ (n-1) $
orthonormal $n$-dimensional vectors
orthogonal to each other and orthogonal to $y_{\nu} $, for $\nu=
1,\ldots, m $.
Then for all $x \in\mathbb{R}^{n} $, we have
\[
\frac{ \Vert x\Vert^2 \operatorname{trace}
(\sigma(x) \sigma^\prime(x))}{ x^\prime\sigma(x) \sigma^\prime(x)
x} -1 \le\frac{(n-1)\varepsilon^{2}+\delta}{1-\delta} < 1.
\]
This is sufficient for the process $X(\cdot)$ to hit the origin in
finite time.

To exclude this situation, we introduce the \textit{effective
dimension} $\mbox{ED}_{\mathcal A}(\cdot)$ of the elliptic
second-order operator $\mathcal A$ defined in (\ref{eq: 6}),
namely
%
%
\begin{equation}\label{eq: effective dim}
\mbox{ED}_{\mathcal A}(x):= \frac{ \Vert x\Vert^2 \operatorname{trace}
(\sigma(x) \sigma^\prime(x) ) }{ x^\prime\sigma(x) \sigma^\prime(x)
x} = \frac{ \Vert x \Vert^2 \operatorname{trace} (A(x)) }{ x^\prime A(x)x}
\end{equation}
for $x \in\mathbb{R}^n \setminus\{0\}$. This function comes from the
theory of the so-called \textit{exterior Dirichlet problem} for
second-order elliptic partial differential equations, pioneered by
Meyers and Serrin \cite{MS60}. These authors showed that
%
%
\begin{equation}\label{eq: ED is bigger 1}
\inf_{x \in\mathbb{R}^{n} \setminus\{0\}}
\mbox{ED}_{\mathcal A}(x)>2
\end{equation}
is a sufficient condition for the existence of solution to an exterior
Dirichlet problem.
In a manner similar to the proof of Proposition~\ref{prop1}, it is
possible to
show that (\ref{eq: ED is bigger 1}) is sufficient for $ \mathbb{P}_{x_{0}}(
X_{1}(t) = \cdots= X_{n}(t) = 0 \mbox{ for some } t \ge0) = 0 $ since
$ \widetilde{R} (\cdot) $ becomes $ \mbox{ED}_{\mathcal A}(\cdot) $
when the matrix $ D $ is replaced by the identity matrix.
[In this
manner, the function $ \widetilde{R} (\cdot) $ of (\ref{eq: h Theta})
is interpreted as a ``local'' version of the effective dimension.]

With $\sigma(\cdot)$ as in (\ref{eq: 3}), the effective
dimension $\mbox{ED}_{\mathcal A}(\cdot)$ satisfies
\begin{eqnarray*}
\mbox{ED}_{\mathcal A}(x) \ge\min_{\nu=1,\ldots, m} \biggl(
\frac{ \Vert x \Vert^2\operatorname{trace} (\sigma_\nu(x) \sigma_\nu
^\prime(x))}{ x^\prime
\sigma_\nu(x) \sigma_\nu^\prime(x) x} \biggr) \ge\min_{\nu
=1,\ldots, m}
\biggl( \frac{\sum_{i=1}^n \lambda_{i \nu}(x)}{\max_{i=1,\ldots,
n} \lambda
_{i \nu}(x)} \biggr)
\end{eqnarray*}
for $x \in\mathbb{R}^n \setminus\{0\}$ where $\{\lambda_{i\nu
}(\cdot),
i=1,\ldots, n\}$ are the eigenvalues of the matrix-valued functions
$\sigma
_\nu(\cdot) \sigma_{\nu}^{\prime}(\cdot)$, for $ \nu= 1,\ldots,
m $.
Thus $\mbox{ED}_{\mathcal A}(\cdot) > 2$ if
\begin{eqnarray*}\label{cond: sufficient 1}
\inf_{x \in\mathbb{R}^n \setminus\{0\}} \min_{\nu=1,\ldots, m}
\biggl(
\frac{\sum_{i=1}^n
\lambda_{i \nu} (x)}{\max_{i=1,\ldots, n} \lambda_{i \nu}(x)}
\biggr) >
2;
\end{eqnarray*}
this can be interpreted as mandating that
the relative size of the maximum eigenvalue is not too large when
compared to all the other eigenvalues.
\end{rem}


%
\begin{rem}
Friedman \cite{MR0346903} established theorems on the nonattainability
of lower-dimensional manifolds by nondegenerate diffusions.
Let $ \mathcal M $ be a closed $ k$-dimensional $C^{2}$-manifold in
$\mathbb{R}^{n} $ with $ k \le n-1 $. At each point $ x \in\mathcal M $, let $
N_{k+i} (x) $ form a set of linearly independent vectors in $\mathbb
{R}^{n} $
which are \textit{normal to} $ \mathcal M $ at $x $. Consider the 
matrix $\alpha(x):= (\alpha_{ij}(x)) $ with
\[
\alpha_{ij}(x) = \langle A(x) N_{k+i} (x), N_{k+j} (x) \rangle;\qquad
1 \le i, j \le n - k, x \in\mathcal M.
\]
Roughly speaking, the strong solution of (\ref{eq: 1}) under a linear
growth condition and a Lipschitz condition on the coefficients
cannot attain $ \mathcal M $ if $ \operatorname{rank}(\alpha(x))$ $
\ge2 $
holds for all $x \in\mathcal M $. The rank indicates how wide the
orthogonal complement of $ \mathcal M $ is. If the rank is large, the
manifold $ \mathcal M $ is \textit{too thin} to be attained. The
fundamental lemma there is based on the solution $u(\cdot)$ of partial
differential inequality ${\mathcal A} u (\cdot) \le\mu u (\cdot)$
for some $\mu\ge0 $, outside but near $ \mathcal M $ with $\lim
_{\mathrm{dist}(x, \mathcal M) \to\infty} u(x) = \infty$ which is different
from our treatment in the previous sections.

Ramasubramanian \cite{R83,Ramasubramanian88} examined the
recurrence and transience of projections of weak solution to (\ref{eq:
1}) for \textit{continuous} diffusion coefficient $\sigma(\cdot) $
showing that any $(n-2)$-dimensional $C^{2}$-manifold is not hit. The
integral test developed there has an integrand similar to the effective
dimension studied in \cite{MS60} as pointed out by M. Cranston in
\textit{Mathematical Reviews.} 

Propositions~\ref{prop1} and~\ref{prop2} are complementary to these
previous general
results since the coefficients here are allowed to be \textit{piecewise
continuous}; however, they depend on the typical geometric
characteristic on the manifold $ \mathcal Z $ we are interested in.
Since the manifold of interest in this work is the zero set $ \mathcal
Z $ of the function $ s(\cdot) $, the projection $ s(X(\cdot)) $ of the
process and the corresponding effective dimensions $ \mbox
{ED}_{\mathcal A}(\cdot) $ and $ \widetilde{R} (\cdot) $ are studied.
\end{rem}
%
%
\begin{rem}\label{rem: vasillios}
As V.~Papathanakos first pointed out,
the conditions (\ref{eq: condition ED and R}), (\ref{eq: condition on R
2}) in Propositions~\ref{prop1} and~\ref{prop2} are
disjoint, and there is a ``gray'' zone of sets of coefficients
which satisfy neither of the conditions. This is because we compare
with Bessel processes, replacing the $n$-dimensional problem by a
solvable \textit{one}-dimensional
problem. In order to look at a finer structure, we discuss a
special case in the next section by reducing it to a \textit
{two}-dimensional problem. This follows a suggestion of A.~Banner.
\end{rem}

\section{A second approach}\label{sec: second approach}


In this section we discuss a class of weak solutions to equation (\ref
{eq: 1}) with the structure (\ref{eq: 3}) which exhibits ``\textit{no
triple collisions}'' using the $n$-dimensional ranked process and the
$(n-1)$-dimensional reflected Brownian motion on polyhedral domains.

\subsection{Ranked process}
Given a vector process $ X(\cdot):= \{ (X_{1}(t), \ldots,
X_{n}(t));$ $0 \le t < \infty\}$, we define the vector $X_{(\cdot)}:=
\{(X_{(1)}(t), \ldots,
X_{(n)}(t)); 0 \le t < \infty\}$ of \textit{ranked processes} ordered
from largest to smallest by
%
%
\begin{equation}\label{ranks}
X_{(k)} (t):= \max_{1 \le i_{1} <\cdots< i_{k} \le n} ( \min
( X_{i_{1}}(t), \ldots, X_{i_{k}}(t) ) );\qquad0 \le t < \infty,
\end{equation}
for $ k = 1,\ldots, n $. If, for every
$j = 1, \ldots, n-2$, the two-dimensional process
%
%
\begin{equation}\label{eq: differences of three}
(Y_{j}(\cdot), Y_{j+1}(\cdot) )^{\prime}:= \bigl(X_{(j)}(\cdot) -
X_{(j+1)}(\cdot), X_{(j+1)}(\cdot) -
X_{(j+2)}(\cdot) \bigr)^\prime
\end{equation}
obtained by looking at the
``gaps'' among the three adjacent ranked processes
$
X_{(j)}(\cdot), X_{(j+1)}(\cdot), X_{(j+2)}(\cdot) $,
never reaches the
corner $(0,0)^\prime$ of $ \mathbb{R}^2 $, almost surely,
then the process $X(\cdot)$ satisfies
%
%
\begin{equation}\label{eq: no triple collision for some t}
\mathbb{P}_{x_{0}} \bigl( X_i (t) = X_j(t) = X_k (t), \mbox{ for
some } (i,j,k),
t > 0 \bigr) = 0
\end{equation}
for $x_{0} \in\mathbb{R}^{n} \setminus\mathcal Z $.
On the other hand, if for some
$j = 1, \ldots, n-2$ the vector of gaps $( X_{(j)} (\cdot) - X_{(j+1)}
(\cdot), X_{(j+1)} (\cdot) - X_{(j+2)} (\cdot) )^\prime$ does
reach the corner $(0,0)^\prime$ of $ \mathbb{R}^2 $ almost surely, then
we have
\begin{eqnarray*} %
\mathbb{P}_{x_{0}}
\bigl( X_i (t) = X_j(t) = X_k (t), \mbox{ for some } (i,j,k), t > 0
\bigr) = 1;\qquad
x_{0} \in\mathbb{R}^{n}.
\end{eqnarray*}

Thus, we are led to study the ranked process $X_{(\cdot)}$ and its
adjacent differences.
In the following we use the parametric result of Varadhan
and Williams \cite{VW85} on Brownian
motion in a two-dimensional wedge with oblique reflection at the
boundary, and the result
of Williams \cite{MR894900} on
Brownian motion with reflection along the faces of a polyhedral domain.

There is a long list of contributions to the study of attainability of
the origin for the Brownian motion with reflection. 
Recently Delarue \cite{MR2453777} considered the hitting time of a
corner by a reflected diffusion in the square. Rogers \cite
{MR1110159,MR1008701} and Burdzy and Marshall \cite{MR1231985} considered
Brownian motion in a half-space with variable angle of reflection. Here
we consider oblique constant reflection on each face of the polyhedral
region. 

\subsection{Reflected Brownian motion}\label{sec: reflected brownian
motion part 1}

Let $ e_{1}, \ldots, e_{n-1} $ be unit vectors
in $\mathbb{R}^{n-1}$, \mbox{$n \ge3 $}, 
and consider the nonnegative orthant
\[
\mathfrak{S}:= \mathbb{R}^{n-1}_{+}
= \Biggl\{ \sum_{k=1}^{n-1} x_{k} e_{k}{}\dvtx{}
x_{1} \ge0, \ldots, x_{n-1} \ge0 \Biggr\},
\]
whose $(n-2)$-dimensional faces
$\mathfrak{F}_{1}, \ldots, \mathfrak{F}_{n-1}$ are given as
\[
\mathfrak{F}_{i}:= \Biggl\{ \sum_{k=1}^{n-1} x_{k} e_{k}{}\dvtx{}
x_{k} \ge0
\mbox{ for } k =1, \ldots, n-1, x_{i} = 0 \Biggr\};\qquad1 \le i
\le n-1.
\]
Let us denote the $(n-3)$-dimensional faces of intersection by
$ \mathfrak{F}_{ij}^o:= \mathfrak{F}_{i} \cap\mathfrak{F}_{j} $ for
$ 1 \le i < j \le n-1 $
and their union by $ \mathfrak{F}^{o}:= \bigcup_{1\le i < j \le n-1}
\mathfrak{F}_{ij}^o $.\vspace*{1pt}

We define the $(n-1)$-dimensional
reflected Brownian motion $
Y (\cdot):= \{(Y_{1}(t),\ldots, Y_{n-1}(t)); t \ge0\}$
on the orthant
$\mathbb{R}^{n-1}_+$ with zero drift, constant $((n-1)\times(n-1))$
constant variance/covariance matrix $ \mathfrak{A}:= \Sigma\Sigma
^\prime$ and
reflection along the faces of the boundary along constant directions by
%
%
\begin{equation}\label{rbm}
\hspace*{10pt} Y(t) = Y(0) + \Sigma B(t) + \mathfrak{R}L(t);\qquad0
\le t <
\infty,
Y(0) \in\mathbb{R}^{n-1}_{+} \setminus\mathfrak{F}^{o}.
\end{equation}
Here, $\{B(t); 0 \le t < \infty\}$ is $(n-1)$-dimensional
standard Brownian motion starting at the origin of $ \mathbb{R}^{n-1}
$. The
$((n-1) \times(n-1))$ reflection matrix $ \mathfrak{R}$ has all its diagonal
elements equal to one, and a spectral radius
strictly smaller than one. Finally, the components of the
$(n-1)$-dimensional process $ L(t):= (L_{1}(t), \ldots, L_{n-1}(t)); 0
\le t < \infty$, are adapted, nondecreasing, continuous and satisfy $
\int_0^\infty Y_i (t)\, d L_i (t) =0 $ [i.e., $ L_i(\cdot)$ is flat off
the set $\{ t\ge0{}\dvtx{} Y_i(t) = 0 \}$] almost surely, for each $i
= 1,\ldots, n-1$.
Note that, if $Y( t)$ lies on $\mathfrak{F}_{ij}^o = \mathfrak{F}_{i}
\cap\mathfrak{F}_{j} $,
then $Y_{i} ( t) = Y_{j} ( t)= 0$ for $1 \le i \neq j \le n-1$.

Harrison and Reiman \cite{HR81} introduced and constructed this process
pathwise through the multi-dimensional Skorohod reflection problem.

\subsubsection{Rotation and rescaling}

Assume that the constant covariance matrix $ \mathfrak{A}= \Sigma
\Sigma' $ is
positive-definite; let $U$ be a unitary matrix whose columns are the
orthonormal eigenvectors of $ \mathfrak{A}$; and let $\mathfrak{L}$
be the
corresponding diagonal matrix of (positive) eigenvalues such that
$ \mathfrak{L}= U^\prime\mathfrak{A}U $. Define
$\widetilde{Y}(\cdot):= \mathfrak{L}^{-1/2} U Y(\cdot) $ and note
that, by this
rotation and rescaling, we obtain
\begin{eqnarray*}\label{eq: definition of Y tilde}
\widetilde{Y}(t) = \widetilde{Y}(0) + \widetilde{B}
(t) + \mathfrak{L}^{-1/2} U \mathfrak{R}L(t);\qquad0 \le t < \infty,
\end{eqnarray*}
from (\ref{rbm}) where $ \widetilde{B}(t):= \mathfrak{L}^{-1/2} U
\Sigma B(t),
0 \le t < \infty$,
is another standard $(n-1)$-dimensional Brownian motion.
We may regard $\widetilde{Y}(\cdot)$ as reflected Brownian
motion in a new state space
$ \widetilde{\mathfrak{S}}:= \mathfrak{L}^{-1/2} U \mathbb
{R}^{n-1}_+ $.
The transformed reflection matrix
$ \widetilde{\mathfrak{R}}:= \mathfrak{L}^{-1/2} U \mathfrak{R}$ can
be written as
%
%
\begin{eqnarray}\label{eq: N and Q}
\hspace*{30pt}\widetilde{\mathfrak{R}} &=&
\mathfrak{L}^{-1/2} U \mathfrak{R}= (\widetilde{\mathfrak{N}} +
\widetilde{\mathfrak{Q}}) \mathfrak{C}
=( \widetilde{\mathfrak{r}}_1, \ldots, \widetilde{\mathfrak{r}}_{n-1} ),\qquad\mbox{where}\nonumber
\\
\mathfrak{C}&:=& \mathfrak{D}^{-1/2},\qquad
\mathfrak{D}:=\operatorname{diag}(\mathfrak{A}),\qquad
\widetilde{\mathfrak{N}}:= \mathfrak{L}^{1/2} U \mathfrak{C}\equiv
(\widetilde{\mathfrak{n}}_1, \ldots, \widetilde{\mathfrak{n}}_{n-1} ),
\\
\widetilde{\mathfrak{Q}} &:=& \mathfrak{L}^{-1/2} U \mathfrak
{R}\mathfrak{C}^{-1} - \widetilde
{\mathfrak{N}}
\equiv( \widetilde{\mathfrak{q}}_1, \ldots, \widetilde{\mathfrak
{q}}_{n-1} ).\nonumber
\end{eqnarray}
Here $ \mathfrak{D} =\operatorname{diag} (\mathfrak{A}) $ is the $((n-1)\times(n-1))$
diagonal matrix with the same diagonal elements as those of
$\mathfrak{A}=\Sigma\Sigma^{\prime}$ (the variances).
The constant vectors
$\widetilde{\mathfrak{r}}_i, \widetilde{\mathfrak{q}}_i, \widetilde
{\mathfrak{n}}_i$,
$i=1,\ldots, n-1$,
are $((n-1) \times1)$ column vectors.

Since $U$ is an orthonormal matrix that rotates the state space
$ \mathfrak{S}= \mathbb{R}^{n-1}_+ $, and
$ \mathfrak{L}^{1/2} $ is a diagonal matrix which changes
the scale in the positive direction,
the new state space $ \widetilde{\mathfrak{S}} $ is an
$(n-1)$-dimensional polyhedron whose $i$th face
$ \widetilde{\mathfrak{F}}_i:= \mathfrak{L}^{-1/2} U \mathfrak{F}_{i} $
has dimension $(n-2)$, for $i = 1, \ldots, n-1$.

Note that
$\operatorname{diag}( \widetilde{\mathfrak{N}}^{\prime} \widetilde
{\mathfrak{Q}}
) = 0 $
and
$\operatorname{diag}( \widetilde{\mathfrak{N}}^{\prime} \widetilde
{\mathfrak{N}})
= I $,
that is, $ \widetilde{\mathfrak{n}}_{i} $ and $ \widetilde{\mathfrak
{q}}_{i} $
are orthogonal and $ \widetilde{\mathfrak{n}}_{i} $ is a unit vector,
that is,
$ \widetilde{\mathfrak{n}}_{i} ^{\prime} \widetilde{\mathfrak{q}}_{i}
= 0 $ and
$ \widetilde{\mathfrak{n}}_{i} ^{\prime} \widetilde{\mathfrak{n}}_{i}
= 1 $ for
$i = 1, \ldots, n-1$.
Also note that
$\widetilde{\mathfrak{n}}_i$ is the inward unit normal
to the $i$th face $ \widetilde{\mathfrak{F}}_{i} $ of the new
state space $
\widetilde{\mathfrak{S}} $
on which the continuous, nondecreasing process $ L_i(\cdot) $
actually increases, for $i = 1, \ldots, n-1$.
The $i$th face $\widetilde{\mathfrak{F}}_{i}$ can be written as
$\{x \in\widetilde{\mathfrak{S}}{}\dvtx{} \widetilde{\mathfrak{n}}_i
^\prime
x = \mathfrak{b}_i\}$ for some
$\mathfrak{b}_i \in\mathbb{R}$, for $i = 1, \ldots, n-1$.

Moreover, the $i$th column $\widetilde{\mathfrak{r}}_{i}$
of the new reflection matrix $\widetilde{\mathfrak{R}}$ is decomposed into
components that
are normal and tangential to $\widetilde{\mathfrak{F}}_{i} $,
that is,
$\widetilde{\mathfrak{r}}_{i} = \mathfrak{C}_{ii}(\widetilde{\mathfrak{n}}_{i} +
\widetilde{\mathfrak{q}}_{i}) $
for $i = 1, \ldots, n-1$ where $\mathfrak{C}_{ii}$ is the
$(i,i)$-element of the
diagonal matrix $\mathfrak{C}$. 
Since the matrix $\mathfrak{L}^{-1/2} U$ of the
transformation is invertible, we obtain\vspace*{-2pt}
%
%
\begin{equation}\label{eq: geometric identity}
\widetilde Y (\cdot) \in\widetilde{\mathfrak{F}}_{ij}^o:= \widetilde{\mathfrak{F}}_{i}
\cap\widetilde{\mathfrak{F}}_{j}
\quad\Longleftrightarrow\quad Y(\cdot) \in\mathfrak{F}_{ij}^o;\qquad
1 \le i < j \le n-1.
\end{equation}
Thus, in order to decide whether the process $Y(\cdot) $ in (\ref{rbm})
attains $ \mathfrak{ F}^o $, it is enough to decide whether the
transformed process $ \widetilde Y(\cdot) $ attains the set $
\widetilde{\mathfrak{F}}^{o}:= \mathfrak{L}^{-1/2} U \mathfrak{F}^{o}
= \bigcup_{1\le i < j
\le
n-1} \widetilde{\mathfrak{F}}_{ij}^o $.

%
\begin{figure}[b]

\includegraphics{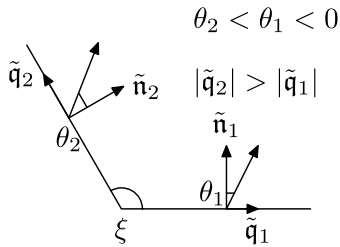}

\caption{Directions of reflection: $\theta_1 + \theta_2 < 0.$}\label
{fig: 1}
\vspace*{-10pt}
\end{figure}

\subsection{Attainability}
With (\ref{eq: geometric identity}) we consider, for $n=3$ and $n > 3$
separately, the hitting times for $ 1 \le i \neq j \le n-1 $:\vspace*{-2pt}
\[
\tau_{ij}:= \inf\{ t > 0{}\dvtx{} Y (t) \in\mathfrak{F}_{ij}^o \}
= \inf\{ t > 0{}\dvtx{} \widetilde Y (t) \in\widetilde{\mathfrak{F}}_{ij}^o
\}.
\]
First we look at the case $n=3$, that is, {\sl two-dimensional
reflected Brownian motion and the hitting time $ \tau_{12} $ of the
origin.} The directions of reflection
$\widetilde{\mathfrak{r}}_{1}$ and $\widetilde{\mathfrak{r}}_{2}$ can
be written in terms
of angles.
Note that the angle $\xi$ of the two-dimensional
wedge $ \widetilde{\mathfrak{S}} $ is positive and
smaller than $\pi$ since all the eigenvalues of $\mathfrak{A}$ are positive.
Let $\theta_1$ and $\theta_2$ with $ -\pi/2< \theta_1, \theta_2 <
\pi/2$
be the angles between $\widetilde{\mathfrak{n}}_1$ and $\widetilde{\mathfrak{r}}_1$, and
between $\widetilde{\mathfrak{n}}_2$ and $\widetilde{\mathfrak{r}}_2 $,
respectively,
measured so that $\theta_1$ is positive if and only
if $ \widetilde{\mathfrak{r}}_1 $ points toward the corner with local
coordinate
$(0,0)^\prime$; similarly for $\theta_2 $.
See Figure~\ref{fig: 1}.

Paraphrasing the result of Varadhan and Williams \cite{VW85}
for Brownian motion reflected on the
\textit{two}-dimensional wedge,
we obtain the following 
result on the relationship between the stopping time and
the sum $ \theta_i + \theta_j $ of angles of reflection directions
when $ n-1=2 $.
%
\begin{lm}[(Theorem~2.2 of \cite{VW85})]\label{lm6}
Suppose that $ \widetilde{Y}(0) = \tilde{y}_{0}
\in\widetilde{\mathfrak{S}}\setminus\widetilde{\mathfrak{F}}^{o} $,
and consider the
ratio $ \beta:= (\theta_{1} + \theta_{2})/\xi$.

The submartingale problem for the reflected Brownian motion on the
two-dimensional wedge is well-posed for $ \beta< 2 $ whereas it has no
solution for $ \beta\ge2 $. If $ 0 < \beta< 2 $, we have $ \mathbb{P}(
\tau_{12} < \infty) = 1 $; if, on the other hand, $ \beta\le0 $,
then we have $ \mathbb{P}( \tau_{12} < \infty) = 0 $.
\end{lm}

In terms of the reflection vectors $ \widetilde{\mathfrak{n}}_{1},
\widetilde{\mathfrak{r}}_{1} $ and
$ \widetilde{\mathfrak{n}}_{2}, \widetilde{\mathfrak{r}}_{2} $, and
with the aid of (\ref{eq: geometric identity}), we can cast this
result as follows; the
proof is in Section~\ref{proof: Proposition 3}.
%
\begin{lm}\label{lm: first lemma of VW}
Suppose that $Y(0) = y_{0} \in\mathbb{R}^{2} \setminus\mathfrak
{F}^{o} $.
If $ \widetilde{\mathfrak{n}}_1^\prime\widetilde{\mathfrak{q}}_2 +
\widetilde{\mathfrak{n}}_2^\prime\widetilde{\mathfrak{q}}_1 > 0 $,
then we have $ \mathbb{P}( \tau_{12} < \infty) = 1 $.
If, on the other hand, $ \widetilde{\mathfrak{n}}_1^\prime\widetilde
{\mathfrak{q}}_2 +
\widetilde{\mathfrak{n}}_2^\prime\widetilde{\mathfrak{q}}_1 \le0 $,
then\break \mbox{$ \mathbb{P}( \tau_{12} < \infty) = 0 $}.
\end{lm}

We consider the general case $ n > 3 $ next.
From (\ref{eq: geometric identity}) and Theorem~1.1 of
Williams \cite{MR894900} we obtain the
following result, valid for $n \ge3$.
%
\begin{lm}\label{lm7}
Suppose that $Y(0) = y_0 \in\mathbb{R}_{+}^{n-1}\setminus\mathfrak
{F}^{o} $
and $n \ge3$ and that the so-called \textup{skew-symmetry condition}
%
%
\begin{equation}\label{SkSym}
\widetilde{\mathfrak{n}}_i^\prime\widetilde{\mathfrak{q}}_j +
\widetilde{\mathfrak{n}}_j^\prime\widetilde{\mathfrak{q}}_i =
0;\qquad
1\le i < j \le n-1,
\end{equation}
holds. Then we have $
\mathbb{P}( \tau< \infty) = 0 $ where $
\tau:= \inf\{ t > 0{}\dvtx{} Y(t) \in\mathfrak{F}^{o} \} $.

Moreover, the components of the adapted, continuous and nondecreasing
process $ L(\cdot) $ defined in $(\ref{rbm})$ are identified then as
the local times at the origin of the one-dimensional component processes
\[
2L_{i}(t) = Y_{i}(t) - Y_{i}(0) - \int^{t}_{0} \operatorname
{sgn}(Y_{i}(s))\, d Y_{i} (s);\qquad
0 \le t < \infty,\mbox{ } i = 1, \ldots, n.
\]
\end{lm}
%
%
\begin{rem}
In the planar (two-dimensional) setting of Lemma~\ref{lm: first
lemma of VW}, the skew-symmetry condition (\ref{SkSym}) takes a weaker
form, that of
an inequality. In the next section we shall discuss some details of the
resulting model as an application of Lemma~\ref{lm7}.
\end{rem}

Lemmata~\ref{lm: first lemma of VW} and~\ref{lm7} lead to
the following result, proved in Section~\ref{subsubsec: application
of Lemma}, on the absence of triple-collisions for a system of $n$
one-dimensional Brownian
particles interacting through their ranks. Let us introduce a
collection $ \{ Q_k^{(i)} \}_{1 \le i,k\le n} $ of polyhedral domains
in $ \mathbb{R}^n $, such that $ \{ Q_k^{(i)} \}_{1 \le i \le n} $ is partition
$ \mathbb{R}^n $ for each fixed~$k$, and $ \{ Q_k^{(i)} \}_{1 \le k\le
n} $ is
partition $ \mathbb{R}^n $ for each fixed $i$. By analogy with (\ref{ranks}),
the interpretation is as follows:
\[
y = (y_1,\ldots, y_n)' \in Q_k^{(i)}\qquad\mbox{means that } y_i
\mbox{ is ranked $k$th among } y_1,\ldots, y_n
\]
with ties resolved by resorting to the smallest index for the highest
rank.
%
\begin{prop}\label{prop3}
For $n \ge3 $, consider the weak solution %
of the equation $(\ref{eq: 5})$ with diffusion coefficient $(\ref{eq:
3})$ where $\sigma(\cdot)$ is the diagonal matrix
%
%
\begin{equation}\label{eq: 70}
\sigma(x):=
\operatorname{diag} \Biggl( \sum_{k=1}^n \widetilde{\sigma}_k
1_{Q_k^{(1)}}(x), \ldots,
\sum_{k=1}^n \widetilde\sigma_k 1_{Q_k^{(n)}}(x) \Biggr);\qquad
x \in\mathbb{R}^{n}.
\end{equation}
If the positive constants
$ \{ \widetilde\sigma_k; 1 \le k \le n \}$
satisfy the linear growth condition
%
%
\begin{equation}\label{growth}
\widetilde{\sigma}_2^2 - \widetilde{\sigma}_{1}^2 =
\widetilde{\sigma}_3^2 - \widetilde{\sigma}_{2}^2 =
\cdots=
\widetilde{\sigma}_{n}^2 - \widetilde{\sigma}_{n-1}^2,
\end{equation}
then $(\ref{eq: no triple collision for some t})$ holds: there are no
triple-collisions among the $n$ 
particles.

If $n = 3$,
the weaker condition $ \widetilde{\sigma}_{2}^{2} - \widetilde
{\sigma
}_{1}^{2} \ge
\widetilde{\sigma}_{3}^{2} - \widetilde{\sigma}_{2}^{2} $ is sufficient
for the absence of triple collisions.
\end{prop}

%
\begin{rem}
The special structure (\ref{eq: 70}) has been studied
in the context of Mathematical Finance. Recent work on interacting
particle systems by Pal and Pitman \cite{PP08} clarifies the long-range
behavior of the spacings between the arranged Brownian particles under
the equal variance condition: $ \widetilde\sigma_{1} = \cdots=
\widetilde\sigma_{n} $; the setting of systems with countably many
particle is also studied there, and related work from Mathematical
Physics on competing tagged particle systems is surveyed.
The ``linear growth'' condition (\ref{growth}) should be seen in the
light of Figure~5.5, page~109 in Fernholz \cite{F02}.
\end{rem}

\section{Application}\label{sec: application}
\subsection{Atlas model for an Equity Market}
Let us recall the Atlas model
%
%
\begin{eqnarray}\label{atlas}
\hspace*{20pt}d X_i(t) &=& \Biggl(\sum_{k=1}^n g_k
1_{Q_k^{(i)}}(X(t)) + \gamma\Biggr)\, d t\nonumber
\\
&&{}+ \sum_{k=1}^n \widetilde\sigma_k 1_{Q_k^{(i)}}(X(t))\,dW_i(t);\qquad
\\
&&{}\hspace*{2pt}\mbox{for } 1 \le i \le n, 0 \le t < \infty,
(X_{1}(0), \ldots, X_{n}(0) )^{\prime} = x_{0} \in\mathbb{R}
^{n},\nonumber
\end{eqnarray}
introduced by Fernholz \cite{F02} and studied by Banner, Fernholz and
Karatzas \cite{BFK05}. Here $ X(\cdot) = (X_1(\cdot),\ldots,
X_n(\cdot
))' $ represents the vector the logarithms of asset capitalizations in
an equity market, and we are using the notation of Proposition~\ref
{prop3}. We
assume that the constants $ \tilde\sigma_k >0 $ and $ g_k $,
$ k = 1,\ldots, n $ satisfy the following conditions which ensure that
$ X(\cdot) $ is ergodic:
\begin{eqnarray*}
g_1&<&0,\qquad g_1 + g_2 <0,\qquad\ldots,\\
g_1 + \cdots+g_{n-1} &<& 0,\qquad
g_1+ \cdots+ g_n = 0.
\end{eqnarray*}
The dynamics of (\ref{atlas}) induce corresponding dynamics for the
ranked processes $ X_{(1)}(\cdot) \ge X_{(2)}(\cdot) \ge\cdots\ge
X_{(n)}(\cdot) $ of (\ref{ranks}).\vspace*{1.5pt} These involve the local times $
\Lambda^{k, \ell} (\cdot)\equiv L^{X_{(k)} - X_{(\ell)}} (\cdot) $ for
$ 1 \le k <\ell\le n $, where 
$ L^Y (\cdot) $ denotes the local time at the origin of a continuous
semimartingale $ Y(\cdot) \ge0 $. An increase in $ \Lambda^{k, \ell}
(\cdot) $ is due to a simultaneous collision of $ \ell- k+1 $
particles in the ranks $ k $ through $ \ell$. In general, when
multiple collisions can occur, there are $(n-1) n/2$ such possible
local times; all these appear then in the dynamics of the ranked
processes, as in Banner and Ghomrasni \cite{BG07}.

Let $ S_{k}(t):= \{ i{}\dvtx{} X_{i} (t) = X_{(k)}(t) \} $ be the set of
indices of processes which are $k$th ranked, and denote its cardinality
by $ N_{k}(t):= \vert S_{k}(t) \vert$ for $0 \le t < \infty$.
Banner and Ghomrasni show in Theorem~2.3 of \cite{BG07} that for any
$n$-dimensional continuous semimartingale $ X(\cdot) = (X_{1}(\cdot),
\ldots, X_{n}(\cdot)) $, its ranked process $ X_{(\cdot)}(\cdot) $ with
components $X_{(k)}(t) = X_{p_{t}(k)}(t)$, $k = 1, \ldots, n$, is
\begin{eqnarray}\label{eq: BG}
d X_{(k)}(t) &=& \sum_{i=1}^{n} \mathbf{1}_{\{X_{(k)}(t) = X_{i}(t)\}
}\, dX_{i}(t)\nonumber
\\[-8pt]\\[-8pt]
&&{}+ \frac{1}{ N_{k}(t) } \Biggl[ \sum_{j=k+1}^{n} d \Lambda^{k, j}(t)
- \sum_{j=1}^{k-1} d \Lambda^{j, k} (t) \Biggr].\nonumber
\end{eqnarray}
Here $ p_{t}:= \{ (p_{t}(1), \ldots, p_{t}(n)) \} $ is the random
permutation of $\{1, \ldots, n\} $ which describes the relation between
the indices of $ X(t) $ and the ranks of $ X_{(\cdot)}(t) $ such that $
p_{t}(k) < p_{t}(k+1) $ if $X_{(k)}(t) = X_{(k+1)}(t) $ for $0 \le t <
\infty$.

Let $ \Pi_n $ be the symmetric group of permutations of $\{ 1, \ldots,
n \} $. The map $ p_{t}{}\dvtx{} \Omega\times[0, \infty) \rightarrow
\Pi_n $ is
measurable with respect to $\sigma$-field generated by the adapted
continuous process $ \{ X(s), 0 \le s \le t\} $, so is predictable.
Consider the inverse map $p_{t}^{-1}:= (p_{t}^{-1}(1), \ldots,
p_{t}^{-1}(n)){}\dvtx{} \Omega\times[0, \infty) \rightarrow\Pi_n
$, also
predictable, indicating the rank of $ X_{i}(t) $ in the $n$-dimensional
vector $X(t) $;
%
%
\begin{equation}\label{eq: p inverse}
X_{(p^{-1}_{t}(i))} (t)= X_{i}(t);\qquad
i = 1, \ldots, n, \mbox{ }0 \le t < \infty.
\end{equation}
%

Under the assumption of ``\textit{no triple collisions}''
[that is, when the only nonzero change-of-rank local times are those of
the form
$\Lambda^{k, k+1} (\cdot) $, $1 \le k \le n-1 $], Fernholz \cite{F02}
considered the stochastic differential equation of the vector of ranked
process $X_{(\cdot)}$ in a general framework; Banner, Fernholz and
Karatzas \cite{BFK05} obtained a rather complete analysis of the Atlas
model (\ref{atlas}).

In this section we apply the main results of the previous
sections to the Atlas model. There
are some cases
of piecewise constant diffusion coefficients which
satisfy the conditions in Proposition~\ref{prop1} or~\ref{prop3}.
Obviously, if the $\{\widetilde\sigma_k^{2}\}$ are all equal,
we are in the case of standard Brownian motion. A bit more
interestingly, if $\{\widetilde\sigma_k^{2} \}$ are linearly growing
in the sense of
(\ref{growth}), we can construct a weak
solution to (\ref{atlas}) with no collision of three or
more particles.
%
\begin{rem}
On page 2305, the paper by Banner, Fernholz and Karatzas \cite{BFK05}
contains the erroneous statement that the ``uniform nondegeneracy of
the variance structure and boundedness of the drift coefficients''
preclude triple collisions. Part of our motivation in undertaking the
present work was a desire to correct this error.
\end{rem}

\subsection{Construction of weak solution}

\subsubsection{Reflected Brownian motion}

Let us start by writing the dynamics of
the sum (total log-capitalization) $ \mathfrak{X}(t):=X_1(\cdot)
+\cdots+X_n(\cdot) $ as
%
%
\begin{equation}\label{eq: sum of X}
\hspace*{30pt}d \mathfrak{X}(t) = n \gamma\, d t + \sum_{i=1}^n
\sum_{k=1}^n \widetilde\sigma_k
1_{Q^{(i)}_k} (X(t))\, d W_i (t) = n \gamma\, d t +
\sum_{k=1}^n \widetilde\sigma_k\, d B_k(t), 
\end{equation}
where
$ B(\cdot):=\{ (B_1 (t),\ldots, B_n(t))', 0\le t < \infty\}$ is given
by $ B_k (t):=\break\sum_{i=1}^n \int_0^t
1_{Q^{(i)}_k} (X(s))\, d W_i (s) $ for $1 \le k \le n $, $ 0 \le t <
\infty$. By the 
F.~Knight theorem (e.g., Chapter~3 in Karatzas and Shreve \cite{KS91}),
this process $ B(\cdot) $ is an $n$-dimensional Brownian motion started
at the origin.

Next, let $h$ and $\widetilde\Sigma$ be the $(n-1)\times1$ vector
and the $(n-1)\times n$ triangular matrix with entries
\begin{eqnarray*}
h:= (g_1-g_2,\ldots, g_{n-1}-g_n)^\prime,\qquad\widetilde
\Sigma:=
\pmatrix{
\widetilde\sigma_{1} & - \widetilde\sigma_{2} & & & \cr
& \widetilde\sigma_{2} & - \widetilde\sigma_{3} & & \cr
& & \ddots& \ddots& \vspace*{2pt}\cr
& & & \widetilde\sigma_{n-1} & - \widetilde\sigma_{n}
},
\end{eqnarray*}
where the elements in the lower-triangular part and
the upper-triangular part, except the first diagonal
above the main diagonal, are zeros.
Then the process $\{ht + \widetilde\Sigma B(t), 0 \le t < \infty\}$
is an $(n-1)$-dimensional Brownian motion starting at the origin of
$ \mathbb{R}^{n-1} $ with constant drift $h$ and the covariance matrix
%
%
\begin{equation}\label{sigma}
\mathfrak{A}:= \widetilde\Sigma\widetilde\Sigma^\prime:=
\pmatrix{
\widetilde\sigma_1^2 + \widetilde\sigma_2^2 & -\widetilde\sigma
_2^2 &&&\cr
-\widetilde\sigma_2^2 & \widetilde\sigma_2^2 + \widetilde\sigma
_3^2&\ddots&\cr
& \ddots& \ddots& -\widetilde\sigma^2_{n-1}\cr
& & -\widetilde\sigma_{n-1}^2&\widetilde\sigma_{n-1}^2+ \widetilde
\sigma_{n}^2
}.
\end{equation}

Now we construct as in Section~\ref{sec: reflected brownian motion part
1} an $(n-1)$-dimensional process
$ Z(\cdot):= \{ (Z_1(t),\ldots, Z_{n-1}(t))^{\prime}, 0 \le t <
\infty\}$ on $\mathbb{R}^{n-1}_{+}$ by
\begin{eqnarray}\label{eq: def of Z reflected BM}
\hspace*{40pt} Z_k(t)&:=& (g_k - g_{k+1})t +
\widetilde\sigma_k B_k(t) - \widetilde\sigma_{k+1} B_{k+1}(t)
\nonumber\\[-8pt]\\[-8pt]
&& {} +\Lambda^{k,k+1}(t) - \tfrac{1}{2}
\bigl( \Lambda^{k-1, k}(t) + \Lambda^{k+1, k+2}(t) \bigr);\qquad
0 \le t < \infty,\nonumber
\end{eqnarray}
for $ k =1,\ldots, n-1 $. Here $ \Lambda^{k,k+1}(\cdot) $ is a
continuous, adapted and nondecreasing process with $ \Lambda
^{k,k+1}(0)=0 $ and $ \int_0^\infty Z_{k}( t)\, d\Lambda^{k,k+1}(t) = 0 $
almost surely. Setting $\Lambda^{0,1}(\cdot) \equiv\Lambda^{n,n+1}
(\cdot) \equiv0$, we write in matrix form
\[
Z(t) = h t + \widetilde\Sigma B(t) + \mathfrak{R}\Lambda(t);\qquad
0 \le t < \infty.
\]
Here $\Lambda(\cdot) =
(\Lambda^{1, 2}(\cdot),\ldots,
\Lambda^{k-1, k}(\cdot))'$ and
the reflection matrix $ \mathfrak{R}=I-\mathfrak{Q}$ is
%
%
\begin{equation}\label{gothRQ}
\mathfrak{R}= I- \mathfrak{Q}:=
\pmatrix{
1 & -1/2 & & &\cr
-1/2 & 1 & \ddots& &\cr
& \ddots& \ddots& \ddots& \cr
& & \ddots& \ddots& -1/2 \vspace*{2pt}\cr
& & & -1/2 & 1
}
.
\end{equation}

\textit{If the process $X(\cdot)$ has no} ``\textit{triple
collisions},'' then
from (\ref{eq: BG}) we get
\begin{eqnarray*}
d X_{(k)} (t) &=& \sum_{i=1}^{n} \mathbf{1}_{\{X_{i} (t) =
X_{(k)}(t)\}}\,dX_{i}(t)\nonumber
\\
&&{}+ \frac{1}{2} \bigl( d \Lambda^{k, k+1} (t)- d \Lambda^{k-1,
k}(t) \bigr),\qquad0 \le t < \infty.\nonumber
\end{eqnarray*}
Substituting (\ref{atlas}) into this equation and subtracting, we
obtain that 
%
%
\begin{equation}\label{eq: the difference of X 1}
X_{(k)}(t) - X_{(k+1)}(t) = Z_{k}(t);\qquad1 \le k \le n-1, 0 \le t <
\infty,
\end{equation}
and that $\Lambda^{k,k+1}(\cdot)$ is the local time at the origin of
the one-dimensional process $Z_k(\cdot) \ge0$ for $k = 1, \ldots,
n-1 $.
In general, the process $X(\cdot) $ may have triple (or
higher-multiplicity) collisions, so that we have additional terms in
(\ref{eq: the difference of X 1}):
%
%
\begin{equation}\label{eq: the difference of X 2}
\hspace*{20pt}X_{(k)}(t) - X_{(k+1)}(t) = Z_{k}(t) + \zeta
_{k}(t),\qquad
1\le k \le n-1, 0 \le t < \infty.
\end{equation}
The contribution $\zeta(\cdot):= (\zeta_{1}(\cdot), \ldots, \zeta
_{n-1}(\cdot)) $ from triple or higher-multiplicity collisions can be written
for $1\le k \le n-1, 0 \le t < \infty$ as $ \zeta_{k} (0)=0 $ and
\begin{eqnarray*}
d \zeta_{k} (t) &=& \sum_{\ell=3}^{n}\ell^{-1} \mathbf{1}_{\{
N_{k}(t) =
\ell\}} \Biggl[ \sum_{j=k+2}^{n} d \Lambda^{k,j} (t)- \sum
_{j=1}^{k-2} d
\Lambda^{j,k}(t) \Biggr] \\
&&{} - \sum_{\ell=3}^{n}\ell^{-1} \mathbf{1}_{\{N_{k}(t) = \ell
\}} \Biggl[ \sum_{j=k+3}^{n} d \Lambda^{k+1,j} (t)- \sum
_{j=1}^{k-1} d \Lambda
^{j,k+1}(t) \Biggr].
\end{eqnarray*}
%

%
\begin{rem}\label{rem: linear local times}
Note that $\zeta(\cdot)
$ consists of (random) linear combinations of local times from
collisions of three or more particles. It is \textit{flat}, unless
there are triple collisions; that is, $\int^{\infty}_{0} \mathbf
{1}_{\mathfrak{ G}^{c}}\, d \zeta(s) = 0 $, 
where $ \mathfrak{ G}:= \{s\ge0{}\dvtx{} X_{i}(t) = X_{j}(t) =
X_{k}(t) \mbox{
for some } 1\le i < j < k \le n \} $. We use this fact with Lemma~\ref
{lm: skorohod type} in the next subsection.
\end{rem}

\subsubsection[Proof of Proposition 3]{Proof of Proposition~\protect\ref{prop3}}\label
{subsubsec: application of Lemma}

Under the assumption of Proposition~\ref{prop3},
we can apply Lemma~\ref{lm7} to obtain
%
%
\begin{equation}\label{eq: cannot reach}
\mathbb{P}\bigl(Z_i (t) = Z_j(t) = 0, \exists t > 0, \exists(i,j),1
\le i
\neq j \le n\bigr) = 0;
\end{equation}
see 
Section~\ref{sec: coefficients structure}.
Thus $ Z(\cdot) $ is a special case of multi-dimensional
reflected Brownian motion for which each continuous, nondecreasing
process $ \Lambda^{k,k+1}(\cdot) $ is exactly the local time at the
origin of 
$ Z_{k}(\cdot) $.

Now let us state the following lemma to examine the local times from
collisions of three or more particles.
Its proof is in Section~\ref{sec: proof of lm: skorohod type}.
%
\begin{lm}\label{lm: skorohod type}
Let $\alpha(\cdot) = \{ \alpha(t);
0 \le t < \infty\} $ be a nonnegative continuous function with
decomposition $ \alpha(t) = \beta(t) + \gamma(t) $ where $\beta
(\cdot)
$ is strictly positive and continuous, and 
$\gamma(\cdot) $ is of finite variation and flat off $ \{ t \ge
0{}\dvtx{}
\alpha(t) = 0 \} $, that is, $\int^{\infty}_{0} \mathbf{1}_{\{
\alpha
(t) > 0\}}\, d \gamma(t) = 0 $. Assume $\gamma(0) = 0 $ and $\alpha
(0) =
\beta(0)>0 $; then, $\gamma(t) = 0 $ and $ \alpha(t) = \beta(t) $ for
all $ 0 \le t < \infty$.
\end{lm}

Under the assumption of Proposition~\ref{prop3}, applying the above
Lemma~\ref
{lm: skorohod type} with (\ref{eq: the difference of X 2}), (\ref{eq:
cannot reach}) and $\alpha(\cdot) = X_{(k)}(\cdot, \omega) -
X_{(k+2)}(\cdot, \omega) $, $\beta(\cdot) = Z_{k}(\cdot, \omega) +
Z_{k+1}(\cdot, \omega) $ and $\gamma(\cdot) = \zeta_{k}(\cdot,
\omega)
+ \zeta_{k+1}(\cdot, \omega) $ for $\omega\in\Omega$, we obtain
$\alpha(\cdot) = \beta(\cdot) $:
%
%
\begin{equation}\label{eq: X k and k+2}
X_{(k)} (\cdot) - X_{(k+2)}(\cdot) = Z_{k}(\cdot) + Z_{k+1}(\cdot
),\qquad
k = 1, \ldots, n-2.
\end{equation}
%
Combining (\ref{eq: X k and k+2}) with (\ref{eq: cannot reach}), we
obtain $X_{(k)}(\cdot) - X_{(k+2)}(\cdot) > 0$ or
\begin{eqnarray*}\label{eq: no triple in rank}
\mathbb{P}\bigl(X_{(k)} (t) = X_{(k+1)}(t) = X_{(k+2)}(t), \exists t > 0,
\exists k,
1 \le k \le n-2 \bigr) = 0.
\end{eqnarray*}
Therefore, there are ``\textit{no triple collisions}'' under the
assumption of Proposition~\ref{prop3}, whose proof is now complete.

\subsubsection{Recovery}\label{subsubsec: recovery}

In conclusion, we \textit{recover} the $n$-dimensional
ranked process $X_{(\cdot)}$
of $X$ by considering a linear transformation.
Specifically, we construct the $n$-dimensional ``ranked'' process,
\[
\Psi_{(\cdot)}(t):= \bigl(\Psi_{(1)}(t),\ldots, \Psi_{(n)}(t)
\bigr);\qquad
0\le t < \infty,
\]
from the sum
$ \mathfrak{X} (t), 0 \le t < \infty$,
defined in (\ref{eq: sum of X}) and
the reflected Brownian motion $Z(\cdot)$, so that
the differences (gaps) satisfy
%
%
\begin{equation}\label{eq: the difference identity}
\Psi_{(k)}(t) - \Psi_{(k+1)}(t) =Z_k(t),\qquad k= 1, \ldots, n-1,
\end{equation}
and
the sum satisfies
%
%
\begin{equation}\label{eq: sum of Psi and X}
\sum_{k=1}^n \Psi_{(k)}(t)= \mathfrak{X} (t);\qquad0 \le t < \infty.
\end{equation}
In particular, each component of $ \Psi_{(\cdot)} (t) $ is uniquely
determined by
\begin{eqnarray*}
\pmatrix{
\Psi_{(1)}( t) \vspace*{2pt}\cr
\Psi_{(2)}( t) \cr
\vdots\cr
\Psi_{(n)}( t)
}
= \frac{1}{n}
\pmatrix{
\mathfrak{X} ( t)+ Z_{n-1}( t)+ (n-2) Z_{n-2}( t) + \cdots+ (n-1) Z_1(t)\vspace*{2pt}\cr
\mathfrak{X} ( t) + Z_{n-1}(\cdot)+ (n-2) Z_{n-2} ( t)+ \cdots- Z_1(t) \cr
\vdots\vspace*{2pt}\cr
\mathfrak{X} ( t)- (n-1) Z_{n-1}(t) - (n-2) Z_{n-2} ( t)- \cdots- Z_1
( t)
}
\end{eqnarray*}
for $0 \le t < \infty$.
Under the assumption of Proposition~\ref{prop3}, we obtain
(\ref{eq: cannot reach}) and hence
with (\ref{eq: the difference identity}) we arrive, in the same way as
discussed
in (\ref{eq: no triple collision for some t}), at
\begin{eqnarray*}\label{eq: no triple collision of Psi}
\mathbb{P}\bigl(\Psi_{(k)} (t) = \Psi_{(k+1)}(t) = \Psi_{(k+2)}(t),
\exists t > 0, 1 \le\exists k \le n-2 \bigr) = 0.
\end{eqnarray*}
Thus, the ranked process $ \{ X_{(\cdot)}(t), 0 \le t < \infty\} $
of the original process $X(\cdot)$ without collision
of three or more particles, and the ranked process
$\Psi_{(\cdot)} (\cdot)$ defined in the above, are equivalent, since
both of them have the same sum
(\ref{eq: sum of Psi and X}) and the same nonnegative difference
processes $Z(\cdot)$
identified in (\ref{eq: the difference of X 1}) and (\ref{eq: the
difference identity}).
We may thus view $\Psi_{(\cdot)} (\cdot)$ as the weak solution
to the SDE for the ranked process $X_{(\cdot)}(\cdot) $.
Finally, we define $\Psi(\cdot):= (\Psi_{1}(\cdot), \ldots, \Psi
_{n}(\cdot)) $ where\vspace*{-1pt}
$\Psi_{i}(\cdot) = \Psi_{(p^{-1}_{t}(i))}(\cdot) $ for $i = 1,
\ldots,
n, $ and $p^{-1}_{t}(i) $ is defined in (\ref{eq: p inverse}). Then,
$\Psi(\cdot)$ is the weak solution of SDE (\ref{atlas}). This
construction of solution leads us to the invariance properties of the
Atlas model given in \cite{BFK05} and \cite{note}.

\begin{appendix}

\section*{Appendix}\label{sec: appendix}
\setcounter{equation}{0}
\subsection{\texorpdfstring{Proof of Lemma~\protect\ref{lm: comparison}}{Proof of Lemma 2.1}}\label{sec:
pf of comparison lem}

From the assumption $ x_{0} \in\mathbb{R}^{n} \setminus\mathcal Z $, where
the zero set $ \mathcal Z $ is defined in (\ref{eq: 9}), it follows
that $ \mathfrak{s}(0) = s(X(\Lambda_{0})) > 0 $ and there exists an
integer $
m_{0} $ such that $ m_{0}^{-1} < \mathfrak{s}(0) < m_{0} $. Recall
that with $
\widetilde{R} (X(\Lambda_{\cdot})) = \mathfrak{d}(\cdot) $ and $
s(X(\Lambda_{\cdot})) = \mathfrak{s}(\cdot) $ we obtained (\ref
{bessel}); namely,
\[
\mathfrak{s}(t) = \mathfrak{s}(0) + \int^{t}_{0} \frac{\mathfrak
{d}(u) - 1}{2 \mathfrak{s}(u)}\, d u + \widetilde B(t);\qquad
0 \le t < \infty.
\]

Let us consider first the case $ \overline\mathfrak
{d}:= \operatorname{essup}\sup_{0\le t <
\infty} \mathfrak{d}(t) < 2 $ for (\ref{eq: lem of hit 2}). Define two
continuous functions $ b_{1}(x):= (\overld-1) / (2 x) $ and $
b_{2}(x):= \overld/ (4 x) $ for $ x \in(0, \infty) $. If $
\overline
\mathfrak{d}< 2 $, then $ b_{1}(\cdot) < b_{2}(\cdot) $ in $ (0,
\infty) $.
For each integer $ m \ge m_{0} $, there exists a nonincreasing
Lipschitz continuous function $ f_{m} (\cdot):= (b_{1}(\cdot) +
b_{2}(\cdot)) / 2 $ with Lipschitz coefficient $ K_{m}:= \max_{x \in
[m^{-1}, m]} \vert b_{2}^{\prime}(x) \vert$, such that $
b_{1}(\cdot) \le f_{m}(\cdot) \le b_{2}(\cdot) $ in $ [m^{-1}, m] $.
%

Define an auxiliary Bessel process $ \mathfrak{r}(\cdot) $ of
dimension $ ( \overld+ 2) / 2$ ($< 2$): 
\[
\mathfrak{r}(t):= \mathfrak{s}(0) + \int^{t}_{0} b_{2}(\mathfrak
{r}(u))\, d u + \widetilde B(t);\qquad
0 \le t < \infty.
\]
%
Consider also the increasing sequence of stopping times
%
%
\begin{equation}\label{taus}
\tau_{m}:= \inf\{t \ge0{}\dvtx{} \max[\mathfrak{s}(t), \mathfrak
{r}(t)] \ge m \mbox{ or }
\min[ \mathfrak{s}(t), \mathfrak{r}(t)] \le m^{-1} \}
\end{equation}
for $ m_{0} \le m < \infty$, and $ \tau:= \inf\{ t \ge0: \mathfrak
{r}(t) =
0 \} $. From the property of the Bessel process with dimension strictly
less than $ 2 $, the process $ \mathfrak{r}(\cdot) $ attains the
origin within
finite time; $ \tau_*:= \lim_{m\to\infty} \tau_{m} \le\tau<
\infty
$ holds a.s.

Now take a strictly decreasing sequence $ \{ a_{n} \}_{n=0}^{\infty}
\subset(0, 1] $ with $ a_{0} = 1 $,\break$ \lim_{n\to\infty} a_{n}
= 0 $
and $ \int_{(a_{n}, a_{n-1})} u^{-2}\, d u = n $ for every $ n \ge1 $.
For each $ n \ge1 $, there exists a continuous function $ \rho
_{n}(\cdot) $ on $ \mathbb{R}$ with support in $ (a_{n}, a_{n-1}) $,
so that
$ 0 \le\rho_{n}(x) \le2(nx^{2})^{-1} $ holds for every $ x > 0 $ and
$ \int_{(a_{n}, a_{n-1})} \rho_{n}(x)\, d x = 1 $.\vspace*{-2pt} Then the function $
\psi_{n}(x):= \int^{\vert x \vert}_{0} ( \int^{y}_{0} \rho
_{n}(u)\, d u
)\, d y $; $ x \in\mathbb{R}$, is even and twice continuous
differentiable with
$ \vert\psi_{n}^{\prime} (x) \vert\le1 $ and $ \lim_{n\to\infty}
\psi_{n}(x) = \vert x \vert$ for $ x \in\mathbb{R}$. Define $
\varphi
_{n}(\cdot):= \psi_{n}(\cdot) \mathbf{1}_{(0, \infty)}(\cdot) $.

By combining the properties of $ \varphi_{n}(\cdot) $, $ b_{1}(\cdot),
b_{2}(\cdot) $ and $ f_{m}(\cdot) $, we see that the difference $
\Delta(\cdot):= \mathfrak{s}(\cdot) - \mathfrak{r}(\cdot) $ is a
continuous process with
\begin{eqnarray*}
\varphi_{n}(\Delta(t))
&\le&\int^{t}_{0} \varphi_{n}^{\prime}(\Delta(u)) \bigl(
b_{1}(\mathfrak{s}(u))
-b_{2}(\mathfrak{r}(u)) \bigr)\, d u
\\
&\le& \int^{t}_{0} \varphi_{n}^{\prime}(\Delta(u)) \bigl(
f_{m}(\mathfrak{s}(u))-f_{m}(\mathfrak{r}(u)) \bigr)\, d u
\\
&\le& K_{m} \int^{t}_{0} \varphi_{n}^{\prime}(\Delta(u)) \bigl
(\mathfrak{s}(u) - \mathfrak{r}(u) \bigr)^{+}\, d u \\
&\le& K_{m} \int^{t}_{0} (\Delta(u))^{+}\, d u;\qquad0 \le t \le
\tau_{m}.
\end{eqnarray*}
Letting $ n \to\infty$ we obtain $ (\Delta(t))^{+} \le K_{m} \int
^{t}_{0}(\Delta(u))^{+}\, d u $ for $ 0 \le t \le\tau_{m} $. From %
the Gronwall inequality and the sample-path continuity of $ \mathfrak
{s}(\cdot),
\mathfrak{r}(\cdot) $ in $ [0, \infty) $, we obtain
$ \Delta(\cdot) = \mathfrak{s}(\cdot) - \mathfrak{r}(\cdot) \le0
$ on $ [0,\tau_{m}]
$ for $ m \ge m_{0} $ and
%
%
\begin{equation}\label{eq: coseq as finite tau tilde}
\hspace*{30pt}\mathfrak{s}(\tau_*) = \lim_{t \to\tau_*} \mathfrak
{s}(t) \le\lim_{t \to\tau_*} \mathfrak{r}
(t) = \mathfrak{r}(\tau_* ) \quad\mbox{and}\quad\max[\mathfrak
{s}(\tau_* ), \mathfrak{r}(\tau_* )
] < \infty,
\end{equation}
almost surely. 
On the other hand, from the definition of $ \{ \tau_{m} \} $ we obtain
$ 0 = \mathfrak{r}(\tau_*)\ge\mathfrak{s}(\tau_*) $, thus $
\mathfrak{s}(\tau_* ) = 0 $ and $
\mathfrak{s}(t) \le\mathfrak{r}(t) $ for $ 0 \le t \le\tau_* $,
a.s., so for $ \overline
\mathfrak{d}= \operatorname{essup} \sup_{0 \le t < \infty}
\mathfrak{d}(t) < 2 $ we conclude
\[
\mathbb{Q}_{x_{0}} \bigl( s(X(t)) = 0 \mbox{ for some } t > 0 \bigr)
= \mathbb{Q}_{x_{0}} \bigl( \mathfrak{s}(t) = 0 \mbox{ for some } t
\ge0
\bigr) = 1.
\]
By the strong Markov property of the process $ X(\cdot) $ under $
\mathbb{Q}
$, we obtain
\[
1 = \mathbb{Q}_{x_{0}} \bigl( s(X(t)) = 0, \mbox{ inf. many } t \ge
0 \bigr)
= \mathbb{Q}_{x_{0}} \bigl( \mathfrak{s}(t) = 0, \mbox{ inf. many }
t \ge
0 \bigr).
\]
This gives (\ref{eq: lem of hit 2}) of Lemma~\ref{lm: comparison}.
Moreover, by the formula of the first hitting-time probability density
function for the Bessel process with dimension $ \overld$ in
Elworthy, Li and Yor \cite{MR1725406} and G{\"{o}}ing--Jaeschke and
Yor \cite
{MR1997032}, we obtain
\begin{eqnarray*}
\mathbb{Q}_{x_{0}} \bigl( \mathfrak{s}(t) = 0, \mbox{ for some } t
\in(0,
T] \bigr)
&\ge& \mathbb{Q}_{x_{0}} \bigl( \mathfrak{r}(t) = 0, \mbox{ for
some } t
\in(0, T] \bigr)
\\
&=&1 - \kappa(T; s(x_{0}), \overld),
\end{eqnarray*}
where the tail probability distribution function $ \kappa(\cdot;
\cdot,
\cdot) $ is defined in (\ref{eq: kappa}).
%

We consider next the case of $\underld:= \operatorname{essinf} \inf_{0 \le t < \infty} \mathfrak{d}(t) \ge2 $. Define $ b_{3}(x):=
(\underld - 1) / (2x) $ and $
b_{4}(x):= \underld / (4 x) $
for $ x \in(0, \infty) $. Following a course similar to the previous
case, using $ b_{3}(\cdot), b_{4}(\cdot) $ and defining a
nonincreasing Lipschitz continuous function $ g_{m}(\cdot):=
(b_{3}(\cdot)+b_{4}(\cdot)) / 2 $ with the Lipschitz coefficient $
L_{m}:= \max_{x \in[m^{-1}, m]} \vert b_{3}^{\prime}(x) \vert$
[rather than using $ b_{1}(\cdot), b_{2}(\cdot) $, $ f_{m}(\cdot) $ and
$ K_{m}$], we obtain the reverse inequality $ \mathfrak{ q} (\cdot)
\le
\mathfrak{s}(\cdot) $ in $ [0, \widetilde{\tau}_{m}] $ a.s. Here $
\mathfrak{
q} (\cdot) $ is the Bessel process in dimension $ (\underld +
2) / 2$ ($ \ge2$); namely 
\[
\mathfrak{ q} (t) = \mathfrak{s}(0) + \int^{t}_{0} b_{4}(\mathfrak{
q} (u))\, d u+ \widetilde B(t);\qquad
0 \le t < \infty,
\]
and the stopping times $ \{ \widetilde{\tau}_m \} $ are defined as in
(\ref{taus}) but with $ \mathfrak{ r} (\cdot) $ replaced by $
\mathfrak
{ q} (\cdot) $.
By a well-known property for Bessel processes of dimension at least $ 2
$, the process $ \mathfrak{ q} (\cdot) $ never attains the origin; that
is, $ \mathfrak{ q} (\cdot) > 0 $ on $ [ 0, \infty) $, a.s.

If $ \widetilde\tau_*:= \lim_{m\to\infty} \widetilde{\tau}_{m} <
\infty$, then by analogy with (\ref{eq: coseq as finite tau tilde}),
we obtain $ \mathfrak{s}(\widetilde\tau_* ) \ge\mathfrak{ q}
(\widetilde\tau
_*) > 0 $ and $ \max[ \mathfrak{s}(\widetilde\tau_*), \mathfrak{ q}
(\widetilde\tau_*) ] < \infty$ a.s., and from the construction of $
\{ \widetilde{\tau}_{m}\} $ a contradiction follows: $ 0 = \mathfrak{s}(
\widetilde\tau_* ) > 0 $. Therefore, $ \mathbb{Q}_{x_{0}} (
\mathfrak
{s}(t) > 0 \mbox
{ for } 0 \le t < \infty) = 1 $. This gives (\ref{eq: lem of hit 1})
of Lemma~\ref{lm: comparison} for $ \underld 
\ge2 $.
%

\subsection{\texorpdfstring{Proof of Propositions~\protect\ref{prop1} and \protect
\ref{prop2}}{Proof of Propositions 1 and 2}}\label{sec: pf of Prop 2}

Proposition~\ref{prop1} and the first half of Proposition~\ref{prop2}
are direct
consequences of Lemma~\ref{lm: comparison} and of the reasoning
developed\vspace*{1pt} in Section~\ref{sec: removal of drift}. Note that $
\langle\widetilde\Theta\rangle(t) \ge c_{0} t, t \ge0$, in this
uniformly nondegenerate case. We obtain (\ref{eq: lem of hit 5}),
because $ \mathbb{Q}_{x_{0}}(s(X(t)) = 0$ for some
\mbox{$t \in[0,T]) \ge\mathbb{Q}_{x_{0}} (\mathfrak{s}(u) = 0$}, for some
$u \in[0, c_{0}T]) $.
Under the original probability measure $ \mathbb{P}_{x_{0}} $, because
of the
drift $\mu(\cdot) $, the process $ s(X(\cdot)) $ is a semimartingale
with the decomposition
\begin{eqnarray*}
d s (X(t)) &=& \biggl( \frac{ (R(x) - 1 ) Q(x)}{2 s(x)} + \frac
{x^{\prime} D
D^{\prime} \mu(x)}{s(x)} \biggr) \bigg\vert_{x = X(t)}\, d t + d
\Theta(t)\\
& =& h (X(t))\, d t + d \Theta(t);\qquad0 \le t < \infty,
\end{eqnarray*}
where $ h(\cdot) $, $ \Theta(\cdot) $ are obtained from $ \widetilde{h}
(\cdot) $, $ \widetilde\Theta(\cdot) $ in (\ref{eq: h Theta}) upon
replacing $ \widetilde{R} (\cdot) $ in (\ref{eq: def of Q}) by $
R(\cdot
) $ in (\ref{eq: modified R}) and $ \widetilde W(\cdot) $ in (\ref{eq:
BM from Grisanov}) by $ W(\cdot) $. The comparison with Bessel
processes is then repeated in a similar manner. When $ \sup_{x \in
\mathbb{R}
^{n} \setminus\mathcal Z} R(x) < 2 $, we get (\ref{eq: triple
collision prob one}) and (\ref{eq: lem of hit 4}).

\subsection{\texorpdfstring{Example in Remark~\protect\ref{rem: example for prop
1}}{Example in Remark 2.3}}\label{sec: proof of example 1}

With some computations we obtain the following simplification of the
effective dimension given in (\ref{eq: effective dim}):
\[
\mbox{ED}_{\mathcal A}(x) = 2 + \frac{
\left[
\matrix{
\Vert x \Vert^{2} - 4 a_{12}(x) \cdot x_{1} x_{2} \mathbf{1}_{\mathcal R_{1+} \cup\mathcal R_{1-}}\vspace*{2pt}\cr
- 4 a_{23}(x) \cdot x_{2} x_{3} \mathbf{1}_{\mathcal R_{2+} \cup\mathcal R_{2-}}\vspace*{2pt}\cr
- 4 a_{31}(x) \cdot x_{3} x_{1} \mathbf{1}_{\mathcal R_{3+} \cup\mathcal R_{3-}}
}
\right] }{x^{\prime} A(x) x}\qquad\mbox{for } x \in\mathbb{R}^{3}
\setminus\{0\}
\]
and
\[
R(x) = 2+ \frac{
\pmatrix{
4 a_{12}(x) \cdot
[(x_{1}-x_{2})^{2}+2(x_{2}-x_{3})(x_{3}-x_{1})]\mathbf{1}_{\mathcal R_{1+} \cup\mathcal R_{1-}}\vspace*{2pt}\cr
{} + 4 a_{23}(x) \cdot
[(x_{2}-x_{3})^{2}+2(x_{3}-x_{1})(x_{1}-x_{2})]\mathbf{1}_{\mathcal R_{2+} \cup\mathcal R_{2-}}\vspace*{2pt}\cr
{} + 4 a_{31}(x) \cdot
[(x_{3}-x_{1})^{2}+2(x_{2}-x_{3})(x_{1}-x_{2})]\mathbf{1}_{\mathcal R_{3+} \cup\mathcal R_{3-}}}
}
{x^{\prime} D D^{\prime} A(x) D D^{\prime} x}
\]
for $x \in\mathbb{R}^{3} \setminus\mathcal Z $ where $ R(\cdot) $
is defined
in (\ref{eq: effective dim}) and $ \mathcal{Z} $ is defined in (\ref
{eq: 9}).
Under the specification (\ref{eq: def of a12}), we verify $\mbox
{ED}(\cdot) > 2$ and $R(\cdot) > 2$, since the denominators of the
fractions on the right-hand sides are positive quadratic forms and
their numerators can be written as
\begin{eqnarray*}
&&\Vert x \Vert^{2} - 4 a_{12} (x) x_{1} x_{2}\\
&&\qquad= (1- 4 a_{12}^{2})x_{2}^{2} + x_{3}^{2} + (x_{1} -
2a_{12}x_{2})^{2} > 0;\qquad x \in\mathcal R_{1+} \cup\mathcal
R_{1-},\\
&& 4 a_{12}(x) [ (x_{1} - x_{2})^{2} + 2(x_{2}-x_{3})(x_{3} -x_{1})]\\
&& \qquad=4 a_{12}(x) {\mathfrak f}_{1}(x) > 0;\qquad x \in\mathcal
R_{1+} \cup\mathcal R_{1-},
\end{eqnarray*}
with similar formulas for $x \in\mathcal R_{i+} \cup\mathcal R_{i-}$,
$i = 2, 3 $.

\subsection{\texorpdfstring{Proof of Lemma~\protect\ref{lm: first lemma of VW}}{Proof of Lemma 3.2}}\label{proof: Proposition 3}

We recall the special geometric
structure of orthogonality $\widetilde{\mathfrak{n}}_i^\prime
\widetilde{\mathfrak{q}}_i= 0 $
and $ \Vert\widetilde{\mathfrak{n}}_{i} \Vert= 1$, and observe that
%
%
\begin{eqnarray}\label{48}
& (\widetilde{\mathfrak{N}}^\prime\widetilde{\mathfrak{Q}}+ \widetilde{\mathfrak{Q}}^\prime
\widetilde{\mathfrak{N}})_{i j}
\matrix{
\ge\cr
<
}
0
\quad\Longleftrightarrow\quad
\widetilde{\mathfrak{n}}_{i}
^\prime\widetilde{\mathfrak{q}}_j +
\widetilde{\mathfrak{n}}_{j} ^\prime\widetilde{\mathfrak{q}}_i
\matrix{
\ge\cr
<
}
0\qquad\forall(i,j).
\end{eqnarray}
Note that if $n = 3 $, that is, $n-1=2$, then
$\widetilde{\mathfrak{n}}_i ^\prime
\widetilde{\mathfrak{q}}_j = \Vert\widehat{\mathfrak{q}}_j
\Vert\operatorname{sgn}( - \theta_j) \sin(\xi)$
for $1\le i \ne j \le2 $
where $\operatorname{sgn}(x):= \mathbf{ 1}_{\{x>0\}} - \mathbf{
1}_{\{x <
0\}} $.
The length $ \Vert\widetilde{\mathfrak{q}}_2 \Vert$
of $\widetilde{\mathfrak{q}}_2$
determines the angle $\theta_2$ and vice versa, that is,
\begin{eqnarray*}
\Vert\widetilde{\mathfrak{q}}_i \Vert
\matrix{
\geq\cr
<
}
\Vert\widetilde{\mathfrak{q}}_j \Vert
\quad\Longleftrightarrow\quad
|\theta_i |
\matrix{
\geq\cr
< }
|\theta_j|.
\end{eqnarray*}
With this observation and $0 < \xi< \pi, \sin( \xi) > 0$, we obtain
\begin{eqnarray*}\label{eq: equivalent relation two}
&& \widetilde{\mathfrak{n}}_i^\prime\widetilde{\mathfrak{q}}_j +
\widetilde{\mathfrak{n}}_j^\prime\widetilde{\mathfrak{q}}_i
= \sin(\xi) \bigl( \Vert\widetilde{\mathfrak{q}}_j \Vert
\operatorname{sgn} ( -\theta
_j) +
\Vert\widetilde{\mathfrak{q}}_i \Vert\operatorname{sgn} (-\theta_i)
\bigr)
\matrix{
\geq\cr
< }
0 \\
&&\qquad
\Longleftrightarrow\quad
\beta= (\theta_i + \theta_j)/\xi
\matrix{
\leq\cr
>
}
0;\qquad1 \le i \neq j \le2.
\end{eqnarray*}
Thus, we apply Lemma~\ref{lm6} and obtain Lemma
\ref{lm: first lemma of VW}.

\subsection{\texorpdfstring{Coefficient structure, and proof of (\protect\ref{eq:
cannot reach})}{Coefficient structure, and proof of (4.10)}}\label{sec: coefficients structure}

Next, we consider the case of linearly growing variance
coefficients defined in (\ref{growth}), and recall
the tri-diagonal matrices $ \mathfrak{A}= \widetilde\Sigma\widetilde
\Sigma^\prime$ as in (\ref{sigma}) and $ \mathfrak{R} $ as in (\ref
{gothRQ}). Consider the $(n-1)$-dimensional
reflected Brownian motion $Y(\cdot)$ defined in (\ref{rbm})
with $\Sigma= \widetilde\Sigma$
and 
$\mathfrak{R}$ as in (\ref{gothRQ}). %
%
Such a
pair $(\widetilde\Sigma, \mathfrak{R})$
satisfies 
%
%
\begin{equation}\label{skewsymidentity2}
(2 \mathfrak{D}- \mathfrak{Q}\mathfrak{D}- \mathfrak{D}\mathfrak
{Q}- 2 \mathfrak{A})_{i j} = 0;\qquad
1 \le i, j \le n-1,
\end{equation}
where $\mathfrak{D}$ is the diagonal matrix with the same diagonal
elements as
$\mathfrak{A}$ %
of (\ref{eq: N and Q}),
and $\mathfrak{Q}$ is the $((n-1)\times(n-1))$ matrix whose
first-diagonal elements above and below the main diagonal
are all $1/2$ and other elements are
zeros as in (\ref{sigma}). In fact,
it suffices to consider $j = i + 1, i = 2, \ldots, n-1$, for which the
equalities (\ref{skewsymidentity2}) are
\begin{eqnarray*}
0 - (\widetilde\sigma_{i}^{2} + \widetilde
\sigma_{i+1}^{2}) - (\widetilde\sigma_{i-1}^{2}+
\widetilde\sigma_{i}^{2}) + 4 \widetilde\sigma_{i}^{2} = 0,
\end{eqnarray*}
or equivalently, (\ref{growth}): $
\widetilde\sigma_{i}^{2} - \widetilde\sigma_{i-1}^{2} =
\widetilde\sigma_{i+1}^{2} - \widetilde\sigma_{i}^{2} $ for $
2 \le i \le n-1 $.
Moreover,\vspace*{1.5pt} the equalities (\ref{skewsymidentity2}) are
equivalent to $
(\widetilde{\mathfrak{N}}^\prime\widetilde{\mathfrak{Q}}+ \widetilde{\mathfrak{Q}}^\prime
\widetilde{\mathfrak{N}}
)_{i j} = 0$ in (\ref{48}).
In fact, from (\ref{eq: N and Q}) with $\mathfrak{D}^{1/2}=\mathfrak
{C}^{-1}$ we compute
\begin{eqnarray*}
\widetilde{\mathfrak{N}}^{\prime} \widetilde{\mathfrak{Q}}&=&
\mathfrak{D}^{-1/2} U^\prime\mathfrak{L}^{1/2} \mathfrak{L}^{-1/2}
U \mathfrak{R}\mathfrak{D}^{1/2} -
\widetilde{\mathfrak{N}}^{\prime} \widetilde{\mathfrak{N}}\\
&=& \mathfrak{D}^{-1/2} (I - \mathfrak{Q}) \mathfrak{D}^{1/2} -
\mathfrak{D}^{-1/2} \mathfrak{A}\mathfrak{D}
^{-1/2}, \\
\widetilde{\mathfrak{N}}^{\prime}\widetilde{\mathfrak{Q}}+ \widetilde{\mathfrak{Q}}^{\prime}
\widetilde{\mathfrak{N}}
&=& 2 I - \mathfrak{D}^{-1/2} \mathfrak{Q}\mathfrak{D}^{1/2}
- \mathfrak{D}^{1/2} \mathfrak{Q}\mathfrak{D}^{-1/2} - 2 \mathfrak{D}^{-1/2} \mathfrak{A}\mathfrak{D}^{-1/2}
\end{eqnarray*}
and multiply both from the left and the right by the diagonal matrix
$\mathfrak{D}^{1/2}$ whose diagonal elements are all positive:
%
%
\begin{equation}\label{eq: NQDA identity}
\mathfrak{D}^{1/2} (\widetilde{\mathfrak{N}}^\prime\widetilde{\mathfrak{Q}}+
\widetilde{\mathfrak{Q}}^\prime\widetilde{\mathfrak{N}}) \mathfrak{D}^{1/2} =
2 \mathfrak{D}- \mathfrak{Q}\mathfrak{D}- \mathfrak{D}\mathfrak{Q}- 2 \mathfrak{A}.
\end{equation}
The equality\vspace*{1.5pt} in the relation (\ref{skewsymidentity2})
is equivalent to the so-called \textit{skew-symmetry condition} $
\widetilde{\mathfrak{N}}^\prime\widetilde{\mathfrak{Q}}+ \widetilde{\mathfrak{Q}}^\prime
\widetilde{\mathfrak{N}}= 0 $ introduced and studied by Harrison and Williams in
\cite
{HW87,MR894900}. It follows from (\ref{48}),
(\ref{skewsymidentity2}) and (\ref{eq: NQDA identity})
that the reflected Brownian motion $ Z(\cdot) $
defined in (\ref{eq: def of Z reflected BM}),
under the assumption of Proposition~\ref{prop3}, is such that
any two dimensional process $(Z_i, Z_j)$
\textit{never} attains the corner
$(0,0)^\prime$ for $1 \le i < j \le n-1 $, that is, (\ref{eq:
cannot reach}) \textit{holds.}
Using this fact, we construct a weak solution to
(\ref{atlas}) from the reflected Brownian motion.
This final step is explained as an application
in the last part of Section~\ref{subsubsec: application of Lemma}.

\subsection{\texorpdfstring{Proof of Lemma~\protect\ref{lm: skorohod type}}{Proof of Lemma 4.1}}\label
{sec: proof
of lm: skorohod type}
We fix an arbitrary $T \in[0, \infty) $. Since $\beta(\cdot) $ is
strictly positive, we cannot have simultaneously $ \alpha(t) = \beta(t)
+ \gamma(t) = 0 $ and $\gamma(t) \ge0 $. The continuous function
$\beta(\cdot) $ attains its minimum on $[0, T] $, so
\begin{eqnarray}\label{eq: the last relation}
\{ t \in[0, T]{}\dvtx{} \alpha(t) = 0 \} &=& \{ t \in[0, T]{}\dvtx
{} \alpha(t) = 0,\gamma(t) < 0 \} \nonumber\\[-8pt]\\[-8pt]
& \subset&\Bigl\{ t \in[0, T]{}\dvtx{} \gamma(t) \le- \min_{0\le
s \le T}
\beta(s) < 0 \Bigr\}.\nonumber
\end{eqnarray}
Let us define $t_{0}:= \inf\{ t \in[0, T]{}\dvtx{} \alpha(t) = 0 \}
$ with
$t_{0} = \infty$ if the set is empty. If $t_{0} = \infty$, then
$\alpha(t) > 0 $ for $0 \le t < \infty$; thus, it follows from the
assumptions $\gamma(0) = 0 $ and $\int^{T}_{0} \mathbf{1}_{\{\alpha(t)
> 0\}}\, d \gamma(t) = 0 $ for $0 \le T < \infty$ that $\gamma(\cdot)
\equiv0 $. On the other hand, if $t_{0} < \infty$, then it follows
from the same argument as in (\ref{eq: the last relation}) that
$\gamma
(t_{0}) < - \min_{0 \le s \le t_{0}} \beta(s) < 0 $. This is
impossible, however, since $\alpha(s) > 0 $ for $0 \le s < t_{0} $ by
the definition of $t_{0}$, and hence the continuous function $\gamma
(\cdot)$ is flat on $[0, t_{0})$, that is, $0 = \gamma(0)=\gamma(t_{0}
- ) = \gamma(t_{0}) $. Thus, $t_{0} = \infty$ and $\gamma(\cdot)
\equiv0 $. Therefore, the conclusions of Lemma~\ref{lm: skorohod
type} hold.

\end{appendix}

\section*{Acknowledgments}
We are grateful to Drs. Robert Fernholz, Adrian Banner and Vasileios
Papathanakos for many discussions on this subject and for providing the
motivation to study this problem. We wish to thank Professors Michel
Emery, Marc Yor and Robert Neel for very helpful discussions, and Dr.
Adrian Banner for his very careful reading of the manuscript. We also
give many thanks to the referees for suggestions and comments that
resulted in major improvements to the manuscript, especially the first
two propositions.

%

\printaddresses

\end{document}